\documentclass[12pt]{amsart}
\usepackage{amssymb}
\usepackage{latexsym}
\usepackage{amsmath}
\usepackage{euscript}
\def\sA{{\mathfrak A}}   \def\sB{{\mathfrak B}}   \def\sC{{\mathfrak C}}
\def\sD{{\mathfrak D}}   \def\sE{{\mathfrak E}}   \def\sF{{\mathfrak F}}
\def\sG{{\mathfrak G}}   \def\sH{{\mathfrak H}}   \def\sI{{\mathfrak I}}
\def\sJ{{\mathfrak J}}   \def\sK{{\mathfrak K}}   \def\sL{{\mathfrak L}}
\def\sM{{\mathfrak M}}   \def\sN{{\mathfrak N}}   \def\sO{{\mathfrak O}}
\def\sP{{\mathfrak P}}   \def\sQ{{\mathfrak Q}}   \def\sR{{\mathfrak R}}
\def\sS{{\mathfrak S}}   \def\sT{{\mathfrak T}}   \def\sU{{\mathfrak U}}
\def\sV{{\mathfrak V}}   \def\sW{{\mathfrak W}}   \def\sX{{\mathfrak X}}
\def\sY{{\mathfrak Y}}   \def\sZ{{\mathfrak Z}}

\def\dA{{\mathbb A}}   \def\dB{{\mathbb B}}   \def\dC{{\mathbb C}}
\def\dD{{\mathbb D}}   \def\dE{{\mathbb E}}   \def\dF{{\mathbb F}}
\def\dG{{\mathbb G}}   \def\dH{{\mathbb H}}   \def\dI{{\mathbb I}}
\def\dJ{{\mathbb J}}   \def\dK{{\mathbb K}}   \def\dL{{\mathbb L}}
\def\dM{{\mathbb M}}   \def\dN{{\mathbb N}}   \def\dO{{\mathbb O}}
\def\dP{{\mathbb P}}   \def\dQ{{\mathbb Q}}   \def\dR{{\mathbb R}}
\def\dS{{\mathbb S}}   \def\dT{{\mathbb T}}   \def\dU{{\mathbb U}}
\def\dV{{\mathbb V}}   \def\dW{{\mathbb W}}   \def\dX{{\mathbb X}}
\def\dY{{\mathbb Y}}   \def\dZ{{\mathbb Z}}

\def\cA{{\mathcal A}}   \def\cB{{\mathcal B}}   \def\cC{{\mathcal C}}
\def\cD{{\mathcal D}}   \def\cE{{\mathcal E}}   \def\cF{{\mathcal F}}
\def\cG{{\mathcal G}}   \def\cH{{\mathcal H}}   \def\cI{{\mathcal I}}
\def\cJ{{\mathcal J}}   \def\cK{{\mathcal K}}   \def\cL{{\mathcal L}}
\def\cM{{\mathcal M}}   \def\cN{{\mathcal N}}   \def\cO{{\mathcal O}}
\def\cP{{\mathcal P}}   \def\cQ{{\mathcal Q}}   \def\cR{{\mathcal R}}
\def\cS{{\mathcal S}}   \def\cT{{\mathcal T}}   \def\cU{{\mathcal U}}
\def\cV{{\mathcal V}}   \def\cW{{\mathcal W}}   \def\cX{{\mathcal X}}
\def\cY{{\mathcal Y}}   \def\cZ{{\mathcal Z}}

\def\bB{{\mathbf B}}
\def\bL{{\mathbf L}}
\def\bQ{{\mathbf Q}}
\def\fG{{\mathsf G}}
\def\fK{{\mathsf K}}
\def\fL{{\mathsf L}}
\def\fM{{\mathsf M}}

\def\d1{{\mathcal D}}

\def\wt#1{{{\widetilde #1} }}
\def\wh#1{{{\widehat #1} }}

\def\bm\chi{\mbox{\boldmath$\chi$}}

\def\RE{{\rm Re\,}}
\def\IM{{\rm Im\,}}
\def\Ext{{\rm Ext\,}}

\def\ran{{\rm ran\,}}
\def\cran{{\rm \overline{ran}\,}}
\def\ker{{\rm ker\,}}
\def\dom{{\rm dom\,}}
\def\mul{{\rm mul\,}}

\def\clos{{\rm clos\,}}
\def\col{{\rm col\,}}
\def\dim{{\rm dim\,}}
\def\codim{{\rm codim\,}}
\def\diag{{\rm diag\,}}

\def\tr{{\rm tr\,}}

 \def\ker{{\xker\,}}
\def\lin{{\rm span\,}}
\def\cspan{{\rm \overline{span}\, }}
\def\sgn{{\rm sgn\,}}
\def\supp{{\rm supp\,}}
\def\cotan{{\rm cotan\,}}
\def\arccot{{\rm arccot\,}}

\def\uphar{{\upharpoonright\,}}
\def\Ext{{\rm Ext\,}}
\def\RE{{\rm Re\,}}
\def\IM{{\rm Im\,}}
\def\wt{\widetilde}
\def\f{\varphi}

\def\ovl{\overline}
\def\d1{{\mathcal D}}
\def\qw{\widetilde Q}

\newtheorem{theorem}{Theorem}[section]
\newtheorem{lemma}[theorem]{Lemma}
\newtheorem{proposition}[theorem]{Proposition}
\newtheorem{corollary}[theorem]{Corollary}
\numberwithin{equation}{section} \theoremstyle{definition}
\newtheorem{definition}[theorem]{Definition}
\newtheorem{remark}[theorem]{Remark}
\newtheorem{example}[theorem]{Example}

\begin{document}

\title[Passive discrete-time systems with a normal main operator]
{Passive systems with a normal main operator \\
and quasi-selfadjoint systems}

\author[Yury Arlinski\u{\i}]{Yu.M. Arlinski\u{\i}}
\address{Department of Mathematical Analysis, \\
East Ukrainian National University, \\
Kvartal Molodyozhny 20-A, \\
Lugansk 91034, \\
Ukraine}
\email{yma@snu.edu.ua}

\author[Seppo Hassi]{S. Hassi}
\address{Department of Mathematics and Statistics, \\
University of Vaasa, \\
P.O. Box 700, \\
65101 Vaasa, \\
Finland}
\email{sha@uwasa.fi}

\author[Henk de Snoo]{H.S.V. de Snoo}
\address{Department of Mathematics and Computing Science,\\
University of Groningen, \\
P.O. Box 800,\\
9700 AV Groningen, \\
Nederland}
\email{desnoo@math.rug.nl}

\thanks{This work was partially supported by the Research Institute
for Technology at the University of Vaasa. The first author was also
supported by the Academy of Finland (projects 117617, 123902) and
the Dutch Organization for Scientific Research NWO (B 61--553)}


\begin{abstract}
Passive systems $\tau=\left\{ T,\sM,\sN,\sH\right\}$ with $\sM$ and
$\sN$ as an input and output space and $\sH$ as a state space are
considered in the case that the main operator on the state space is
normal. Basic properties are given and a general unitary similarity
result involving some spectral theoretic conditions on the main
operator is established. A passive system $\tau$ with $\sM=\sN$ is
said to be quasi-selfadjoint if ${\rm ran\,}(T-T^*)\subset\mathfrak
N$. The subclass ${\bf S}^{qs}(\mathfrak N)$ of the Schur class
${\bf S}(\mathfrak N)$ is the class formed by all transfer functions
of quasi-selfadjoint passive systems. The subclass ${\bf
S}^{qs}(\mathfrak N)$ is characterized and minimal passive
quasi-selfadjoint realizations are studied. The connection between
the transfer function belonging to the subclass ${\bf
S}^{qs}(\mathfrak N)$ and the $Q$-function of $T$ is given.
\end{abstract}

\maketitle

\section{Introduction}

Let $\sM,\sN$, and $\sH$ be separable Hilbert spaces and let
\begin{equation}\label{abcd0}
T=\begin{pmatrix} D&C\cr
B&A\end{pmatrix}:\begin{pmatrix}\sM\\
\sH\end{pmatrix}\to
\begin{pmatrix}\sN\\ \sH\end{pmatrix}
\end{equation}
be a bounded linear operator. Here and in the following, it will be
tacitly assumed that the spaces in the righthand side are orthogonal
sums: $\sM \oplus \sH$ and $\sN \oplus \sH$. The system of equations
\begin{equation}
\label{passive} \left\{
\begin{array}{l}
 h_{k+1}=Ah_k+B\xi_k,\\
 \sigma_k=Ch_k+D\xi_k,
\end{array}
\right. \qquad k\ge 0,
\end{equation}
describes the evolution of a \textit{linear discrete time-invariant
system} $\tau=\left\{T,\sM,\sN,\sH\right\}$. The Hilbert spaces
$\sM$ and $\sN$ are called the input and the output spaces,
respectively, and the Hilbert space $\sH$ is called the state space.
The operators $A$, $B$, $C$, and $D$ are called the \textit{main
operator}, the \textit{control operator}, the \textit{observation
operator}, and the \textit{feedthrough operator} of $\tau$,
respectively. The subspaces
\begin{equation}
\label{CO} \sH^c =\cspan\{\,A^{n}B\sM:\,n \in \dN_0 \}
 \mbox{ and } \sH^o =\cspan\{\,A^{*n}C^*\sN:\,n \in \dN_0 \}
\end{equation}
are called the \textit{controllable} and \textit{observable}
subspaces of $\tau=\left\{T,\sM,\sN,\sH\right\}$,
respectively. If $\sH^c=\sH$ ($\sH^o=\sH$) then the system $\tau$ is
said to be \textit{controllable} (\textit{observable}), and
\textit{minimal} if $\tau$ is both controllable and observable. If
$\sH=\clos \{\sH^c+\sH^o\}$ then the system $\tau$ is said to be a
\textit{simple}.  Two discrete-time systems
$\tau_1=\left\{T_1,\sM,\sN,\sH_{1}\right\}$
 and
$\tau_2=\left\{T_2,\sM,\sN,\sH_{2}\right\}$ are  \textit{unitarily
similar} if %$D_1=D_2$ and
there exists a unitary operator $U$ from
$\sH_{1}$ onto $\sH_{2}$ such that
%\[
%A_1 =U^{-1}A_2U,\quad B_1=U^{-1}B_2,\quad C_1=C_2 U.
%\]
\begin{equation}\label{unisim}
 A_2=UA_1U^*,\quad  B_2=UB_1, \quad  C_2=C_1U^*, \mbox{ and } D_2=D_1.
\end{equation}
If the linear operator $T$ is contractive (isometric, co-isometric,
unitary), then the corresponding discrete-time system is said to be
\textit{passive} (\textit{isometric, co-isometric, conservative}).
The \textit{transfer function}
\begin{equation}
\label{TrFu} \Theta (\lambda):=D+\lambda C(I -\lambda A)^{-1}B,
\quad \lambda \in \dD,
\end{equation}
of the passive system $\tau$ in \eqref{passive} belongs to the
\textit{Schur class} ${\bf S}(\sM,\sN)$, i.e., $\Theta (\lambda)$ is
holomorphic in the unit disk $\dD=\{\lambda \in \dC:|\lambda|<1\}$
and its values are contractive linear operators from $\sM$ into
$\sN$.
% As is well known \cite{BrR2}, \cite{SF}, \cite{Ando} \cite{A},
% \cite{ArKaaP} e
Every operator-valued function $\Theta(\lambda)$ from the Schur
class ${\bf S}(\sM,\sN)$  can be realized as the transfer function
of a passive system, which can be chosen as   observable
co-isometric (controllable isometric,  simple conservative, passive
minimal). Moreover two isometric and controllable (co-isometric and
observable, simple conservative) systems having the same transfer
function are unitarily similar.  D.Z.~Arov~\cite{A} has shown that two
minimal passive systems $\tau_1$ and $\tau_2$ with the same transfer
function $\Theta(\lambda)$ are only \textit{weakly similar}, i.e.,
there is a closed densely defined operator $Z:\sH_{1}\to\sH_{2}$
such that $Z$ is invertible, $Z^{-1}$ is densely defined, and
\begin{equation}\label{pseudosim}
Z A_1f =A_2 Zf, \quad C_1f=C_2 Zf,\quad f\in\dom Z,\quad\mbox{and}
\quad ZB_1=B_2.
\end{equation}
Weak similarity preserves neither the dynamical properties of the
system nor the spectral properties of its main operator $A$. In
\cite{ArNu1}, \cite{ArNu2} necessary and sufficient conditions have
been established for minimal passive systems with the same
transfer function to be (unitarily) similar.  In \cite{AHS2} a
parametrization of the contractive block-operator matrices in
\eqref{abcd0} was used to establish some new aspects and some
explicit formulas for the interplay between the system $\tau$, its
transfer function $\Theta (\lambda)$, and the Sz.-Nagy--Foia\c{s}
characteristic function of the contraction $A$.

In this paper the same approach is applied to study
passive systems with a normal main operator, including the class of
\textit{passive quasi-selfadjoint} systems ($pqs$-systems for short),
as defined in the paper.
Furthermore, using the famous Mergelyan's theorem from complex analysis
a general unitary similarity result is proved for such systems.

The
passive system $\tau=\{T,\sN,\sN,\sH\}$ is called a $pqs$-system if
the operator
\begin{equation}\label{abcd}
T=\begin{pmatrix} D&C\cr
B&A\end{pmatrix}:\begin{pmatrix}\sN\\
\sH\end{pmatrix}\to
\begin{pmatrix}\sN\\ \sH\end{pmatrix}
\end{equation}
is a \textit{quasi-selfadjoint contraction} ($qsc$-operator for
short), i.e., $T$ is a contraction and $\ran (T-T^*)\subset \sN$,
cf. \cite{AHS1}. This last condition is equivalent to $A=A^*$ and
$C=B^*$. If $\tau$ is a $pqs$-system, then the transfer function
\eqref{TrFu} of $\tau$ takes the form
\[
\Theta (\lambda)=W(\lambda)+D,
\]
where the function $W(\lambda)$ is a Herglotz-Nevanlinna function
defined on %the domain
$\Ext\{(-\infty,-1]\cup [1,\infty)\}$.
%($\IM \lambda\, (W(\lambda)-W^*(\lambda)\ge 0).$
The subclass ${\bf S}^{qs}(\sN)$ of the Schur class ${\bf S}(\sN)$
of $\bL(\sN)$-valued functions  is the class of all transfer
functions of $pqs$-systems $\tau=\{T;\sN,\sN,\sH\}$.
A necessary and sufficient condition for the function $\Theta(\lambda)$ to be in the class
${\bf S}^{qs}$ is given, the minimal $pqs$-systems with the given
operator-valued function $\Theta(\lambda)$ from the class ${\bf
S}^{qs}$ is constructed using operator representations of Herglotz-Nevanlinna
functions. Moreover, a necessary and sufficient condition for the function
$\Theta(\lambda)\in{\bf S}^{qs}$ to be inner (co-inner) is proved and
connections with $pqs$-system and other minimal systems with the same transfer
function are established. Also it is shown that if, for instance,
$\Theta(\lambda)\in {\bf S}^{qs}(\sN)$ and $\f_\Theta(\lambda)=0$
($\psi_\Theta(\lambda)=0$) then $\Theta(\lambda)$ is inner
(co-inner). A matrix form of the inner function from the
class ${\bf S}^{qs}(\sN)$ when $\dim\sN<\infty$ is also given, and in the case
of scalar functions from the class ${\bf S}^{qs}(\sN)$ a minimal representation
is obtained by means of Jacobi matrices.

\section{Preliminaries}

Let $\sM$ and $\sN$ be Hilbert spaces and let $\Theta(\lambda)$
belong to   the Schur class ${\bf S}(\sM,\sN)$.  The notation
$\Theta(\xi)$, $\xi\in\dT$, stands for the non-tangential strong
limit value of $\Theta(\lambda)$ which exist almost everywhere on
$\dT$, cf. \cite{SF}.  A function $\Theta(\lambda) \in {\bf
S}(\sM,\sN)$ is said to be \textit{inner} if
$\Theta^*(\xi)\Theta(\xi)=I_\sM$ for almost all $\xi\in\dT$, and it
is said to be \textit{co-inner} if $\Theta(\xi)\Theta^*(\xi)=I_\sN$
for almost all $\xi\in\dT$. A function $\Theta(\lambda)\in {\bf
S}(\sM,\sN)$ is said to be \textit{bi-inner} if it is both inner and
co-inner.

\subsection{Contractions and their defect operators}

Let $A\in{\bL}(\sH_1,\sH_2)$ be a contraction, in other words, let
$\|Af\| \le \|f\|$ for all $f \in \sH$, or equivalently $I-A^*A \ge
0$.  The selfadjoint operator $D_A =(I - A^*A)^{1/2}$ is said to be
the \textit{defect operator} of $A$. Observe that $\ker D_A=\ker
D_A^2=\ker (I-A^*A)$, and that
\begin{equation}\label{ker}
 \ker (I-A^*A)= \{\,f \in \sH :\, \|Af\|=\|f\|\,\}.
\end{equation}
Clearly, any contraction $A$ satisfies
% \begin{equation}\label{ker1}
%  \ker (I-A^*A) \subset \ker (I-(A^*A)^n),  \quad n\in\dN.
%  %\quad \ker (I-AA^*)
%  %\subset \ker (I-(AA^*)^n).
% \end{equation}
% Furthermore,  for a contraction one has
\[
 I\ge A^*A \ge\cdots \ge A^{*n}A^n \ge 0,
%\quad I\ge AA^* \ge\cdots \ge A^nA^{*n} \ge 0,
 \quad n\in\dN,
\]
in other words the sequence  $\|A^nf\|$
%and $\|A^{*n}f\|$
with $f \in \sH$ is monotonically nonincreasing. In particular, the
strong limit
\begin{equation}\label{sa}
S_A=s-\lim A^{*n}A^n, %\quad S_{A^*}=s-\lim A^{n}A^{*n}
\end{equation}
exists as an operator in ${\bL}(\sH_1)$, cf. \cite[p. 261]{SF0}.
The defect operators $D_A$ and $D_{A^*}$ satisfy the following
commutation relation:
\begin{equation} \label{comm}
AD_A = D_{A^*}A, \quad D_AA^*=A^*D_{A^*}.
\end{equation}
Let ${\sD}_A$ stand for the closure of the range ${\ran}D_A$.  Then
%and therefore the block operator
\begin{equation}
\label{blockmatrix} %T=
\begin{pmatrix}-A&D_{A^*}\cr
D_{A}&A^*\end{pmatrix}:\begin{pmatrix}\sD_A\\
\sH_2\end{pmatrix}\to
\begin{pmatrix}\sD_{A^*}\\\sH_1\end{pmatrix}
\end{equation}
is unitary. Define the subspaces $\sH_{A,0}$ and $\sH_{A,1}$ by
\[
 \sH_{A,1}=\left\{\, f\in\sH:\, \|f\|=\|A^nf\|=\|A^{*n}f\|,\;
n \in \dN\,\right\}, \quad \sH_{A,0}=\sH \ominus \sH_{A,1}.
\]
Then $\sH=\sH_{A,0}\oplus \sH_{A,1}$ is a canonical orthogonal
decomposition of $\sH$ such that
\begin{equation}\label{Adecom}
A=A_0\oplus A_1, \quad A_j=A\uphar \sH_{A,j}, \quad j=0,1,
\end{equation}
where $\sH_{A,0}$ and $\sH_{A,1}$ reduce $A$, $A_0$ is a completely
non-unitary contraction, and $A_1$ is a unitary operator. The
function
\begin{equation}
\label{CHARFUNC1}
 \Phi_{A}(\lambda):=\left(-A+\lambda D_{A^*}(I_\sH-\lambda
A^*)^{-1}D_{A}\right)\uphar\sD_{A}, \quad \lambda \in \dD,
\end{equation}
is the Sz.-Nagy--Foia\c{s} characteristic function of the
contraction $A$. It belongs to the Schur class ${\bf
S}(\sD_A,\sD_{A^*})$; cf. \cite{SF}. In fact, a straightforward
calculation using the identities $(I-\lambda
A)^{-1}=I+\lambda(I-\lambda A)^{-1}A$ and \eqref{comm} yields
\begin{equation}
\label{chardef}
 D_{\Phi_{A}(\lambda)}^2=(1-\lambda\ovl\lambda)D_A(I-\ovl\lambda
 A)^{-1}(I-\lambda A^*)^{-1}D_A\uphar\sD_{A},
%% \quad \lambda \in \dD: Better to use resolvent set of T here!
\end{equation}
which shows that $\Phi_{A}(\lambda)$ is contractive for $\lambda \in
\dD$. Note also that $\Phi_{A^*}(\lambda)$ is the transfer function
of the conservative system
\begin{equation}
\label{SIGMA} \Sigma=\{T,\sD_{A^*},\sD_A,\sH\},
\end{equation}
where
\begin{equation}
\label{blockmatrix1}
 T=\begin{pmatrix}-A^* & D_{A}\cr
    D_{A^*} & A\end{pmatrix}:\begin{pmatrix}\sD_{A^*}\\
\sH\end{pmatrix}\to
\begin{pmatrix}\sD_{A}\\ \sH\end{pmatrix}.
\end{equation}
Let $A \in \bL(\sH_1,\sH_2)$ be a contraction and let
$\Sigma$
%=\{T,\sD_{A^*},\sD_A,\sH\}$
be the corresponding
conservative system in \eqref{SIGMA}, \eqref{blockmatrix1}. Then
the controllable and observable subspaces, as defined in \eqref{CO},
are given by
\begin{equation}\label{coSig}
\sH^c_\Sigma=\cspan\{A^{n}D_{A^*}\sD_{A^*}:\,n \in \dN_0  \},\quad
\sH^o_\Sigma=\cspan\{A^{*n}D_{A}\sD_{A}:\,n \in \dN_0 \,\}.
\end{equation}
Observe that $A$ is completely nonunitary if and only if $\Sigma$ is
minimal.
Since clearly $\Phi_{A}(\lambda)^*=\Phi_{A^*}(\ovl\lambda)$, one has
also
\begin{equation}
\label{chardef*}
 D_{\Phi_{A}(\lambda)^*}^2=(1-\ovl\lambda\lambda)D_{A^*}(I-\lambda
 A^*)^{-1}(I-\ovl\lambda A)^{-1}D_{A^*}\uphar\sD_{A^*}.
%% \quad \lambda \in \dD: Better to use resolvent set of T here!
\end{equation}
Observe that if $\xi\in\dT:=\{\xi \in \dC: |\xi|=1\}$ belongs to the
resolvent set of $A$, then \eqref{chardef} and \eqref{chardef*} show
that $\Phi_{A}(\xi)$ is a unitary operator; cf. \cite[p. 239]{SF}.

 A contraction $A$ in a Hilbert space
$\sH$ is said to belong to the \textit{classes} $C_{0\,\cdot}$ or
$C_{\cdot\, 0}$ if
\[
s-\lim\limits_{n\to\infty}A^n=0  \quad \mbox{or} \quad
s-\lim\limits_{n\to\infty}A^{*n}=0,
\]
respectively. By definition, $C_{00}:=C_{0\,\cdot}\cap C_{\cdot\,
0}$. Hence $A\in C_{00}$ precisely when
\begin{equation}\label{lim0}
 s-\lim\limits_{n\to\infty}A^n=s-\lim\limits_{n\to\infty}A^{*n}=0.
\end{equation}
Observe that $A \in C_{00}$  implies that $A$ is completely
nonunitary, cf. \eqref{Adecom}. The completely non-unitary part of a
contraction $A$ belongs to the class $C_{\cdot\,0}$, $C_{0\,\cdot}$,
or $C_{00}$ if and only if its characteristic
function $\Phi_A(\lambda)$ in \eqref{CHARFUNC1} is inner, co-inner,
or bi-inner, respectively; cf. \cite[Theorem VI.2.3]{SF}.
 It follows from \eqref{sa} that $S_A=0$ implies
$A\in C_{0 \ \cdot}$ and that $S_{A^*}=0$ implies  $A\in C_{\cdot
\,0}$.

A contraction $A$ is said to be \textit{strict} if $\|Af\|<\|f\|$
for all nontrivial $f\in \sH_1$. Note that in view of \eqref{ker} a
contraction $A$ is strict if and only if $\ker D_A= \ker D_A^2=\ker
(I-A^*A)=\{0\}$. Finally, a passive system $\tau=\{T;\sM,\sN,\sH\}$
is said to be \textit{strongly stable} or \textit{strongly
co-stable} if the main operator $A$ belongs to the class
$C_{0\,\cdot}$ or $C_{\cdot\, 0}$, respectively; see \cite{Arov},
\cite{ArSt}.

\subsection{Some properties of normal contractions}

An operator $A\in{\bL}(\sH_1,\sH_2)$ is said to be \textit{normal}
if $A^*A=AA^*$, or equivalently, if $\|Af\|=\|A^*f\|$ for all $f \in
\sH$, cf. \cite[p. 281]{SF0}.
It is clear from  $\sH=\sH_{A,0} \oplus \sH_{A,1}$
and the orthogonal decomposition in \eqref{Adecom}
that a contraction
$A$ is normal if and only if its completely nonunitary part $A_0$ is normal
in $\sH_{A,0}$.
If $A$ is a normal contraction  then,
parallel to \eqref{ker}, one has
\begin{equation}\label{ker2}
\begin{split}
 &\ker (I-(A^*A)^n)=\ker (I-A^{*n}A^n) =\{\,f \in \sH :\, \|A^nf\|=\|f\|\,\} \\
&=\ker (I-(AA^*)^n) =\ker (I-A^nA^{*n}) =\{\,f \in \sH :\,
\|A^{*n}f\|=\|f\|\,\}.
\end{split}
\end{equation}
Moreover, if $A$ is a normal contraction, then the defect operators
$D_A$ and $D_{A^*}$ satisfy $D_A=D_{A^*}$ and ${\sD}_A={\sD}_{A^*}$;
in addition, \eqref{comm} reads as
\begin{equation} \label{comm+}
 AD_A = D_{A}A, \quad A^*D_A=D_A A^*.
\end{equation}
%In this case many properties of the
%characteristic function of $A$ simplify.

\begin{lemma}\label{NEW}
Let $A\in{\bL}(\sH_1,\sH_2)$ be a normal contraction. Then the
strong limit $S_A$ satisfies $S_A=S_{A^*}$ and
\begin{equation}\label{sa1}
 S_A(I-A^*A)=0.
\end{equation}
If, in addition, $A$ is strict, then $S_A=0$.
\end{lemma}

\begin{proof}
%Clearly, if bounded operators $B$ and $C$ satisfy the inequality $B
%\le C$, then $A^*BA \le A^*CA$  for any bounded operator $A$. Hence,
%for a contraction $A$ it follows that
%\[
% I\ge A^*A \ge\cdots \ge A^{*n}A^n \ge 0, \quad n\in\dN.
%\]
%Hence the limit $S_A$ exists, cf. \cite{NF0}.
If $A$ is normal, then \eqref{sa} implies that $S_A=S_{A^*}$ and
\[
 S_AA^*A= (s-\lim A^{*n}A^{n})A^*A=s-\lim A^{*(n+1)}A^{n+1}=S_A,
\]
which leads to \eqref{sa1}.
\end{proof}

%In the case that the main operator $A$ is normal, many properties of
%the characteristic function of $A$, as well as of the system $\tau$
%and of its transfer function $\Theta(\lambda)$ in \eqref{TrFu},
%simplify.

\begin{proposition}\label{nstable}
Let   $A \in \bL(\sH_1,\sH_2)$ be a normal contraction. Then the
following statements are equivalent:
\begin{enumerate}\def\labelenumi{\rm (\roman{enumi})}
\item $A\in C_{00}$;
\item $A$ is completely non-unitary;
\item $A$ is strict.
\end{enumerate}
Moreover, the characteristic function $\Phi_A(\lambda)$ of $A$ in
\eqref{CHARFUNC1} is bi-inner.
\end{proposition}

\begin{proof}
(i) $\Rightarrow$ (ii) This implication is a general fact for not
necessarily normal contractions.

(ii) $\Rightarrow$ (iii) Let $A$ be completely non-unitary. Assume
that $A$ is not strict. Then there exists an element $0\neq f_0\in
\sH_1$ such that $\|Af_0\|=\|f_0\|$.
Since $\ker (I-A^*A) \subset \ker (I-(A^*A)^n)$, $n\in\dN$,
it follows from \eqref{ker}
%\eqref{ker1},
and \eqref{ker2} that
%Hence, $f_0\in\ker D_A=\ker D_A^2$ and consequently
%$(A^*A)^nf_0=f_0$ for all $n\in\dN$. Since $A$ is normal, this
%implies that
$\|f_0\|=\|A^nf_0\|=\|A^{*n}f_0\|>0$. This contradicts the fact that $A$ is
completely nonunitary.
%Hence $A$ is not a , so that \eqref{lim0} does not hold.

(iii) $\Rightarrow$ (i) Let $A$ be strict, so that $\ker (I-A^*A)=\{0\}$.
Then Lemma \ref{NEW} implies that $S_A=0$, which leads to
\eqref{lim0}, so that $A\in C_{00}$.

Observe that if $A$ is normal then $\sH_{A,1}=\ker D_A$ as was just
shown above. The completely non-unitary part
$A_0$ of $A$ is normal and satisfies $\ker D_{A_0}=\{0\}$. Thus
$A_0\in C_{00}$, i.e., $\Phi_A(\lambda)$ is bi-inner; cf.
\cite[Theorem VI.2.3]{SF}.
\end{proof}

%$f_0\in\ker D_A$ with $\|A^{*n}f_0\|=\|A^nf_0\|=\|f_0\|>0$ for all
%$n\in\dN$; see the proof of Lemma~\ref{nstable}. This means that
%$f_0\in\sH_{A,1}$, so that $A$ is not a completely non-unitary
%contraction; cf. \eqref{Adecom}. This proves the reverse implication
%(v) $\Longrightarrow$ (vi).
%
%
%%Since $I\ge A^*A
%%\ge\cdots \ge A^{*n}A^n \ge 0$ for all $n\in\dN$ the strong limit
%%$S=s-\lim A^{*n}A^n \ge 0$ exists. Moreover, since $A$ is normal,
%%\[
%% SA^*A= (s-\lim A^{*n}A^{n})A^*A=s-\lim A^{*(n+1)}A^{n+1}=S,
%%\]
%%so that $S(I-A^*A)=0$. Hence, if $A$ is a strict contraction, so
%%that $\ker D^2_A=\{0\}$, then $Sf=0$ for all $f\in \sD_A=\sH_1$ and
%%\eqref{lim0} holds.
%
%(ii) $\Rightarrow$ (iii) This is clear from \eqref{Adecom}.
%
%
%
%
%Now assume that \eqref{lim0} holds. If $A$ is not strict, then there
%exists an element

If a contraction $A$ is normal, then its controllable and observable subspaces
coincide, which leads to the following observation.

\begin{proposition}\label{nnstable}
Let $A \in \bL(\sH_1,\sH_2)$ be a normal contraction, and let
$\Sigma$
%=\{T,\sD_{A^*},\sD_A,\sH\}$
be the corresponding
conservative system in \eqref{SIGMA}, \eqref{blockmatrix1}.   Then $
\sH^c_\Sigma=\sH^o_\Sigma$
%\begin{equation}\label{IDE}
% \sH^c_\Sigma=\sH^o_\Sigma,
%\end{equation}
and the following statements are equivalent:
\begin{enumerate}
\def\labelenumi{\rm (\roman{enumi})}
\item $\Sigma$ is simple;
\item $\Sigma$ is controllable;
\item $\Sigma$ is observable;
\item $\Sigma$ is minimal.
%\item $A$ is completely non-unitary.
%\item $A$ is strict. %,\ker D_A=\{0\}$ or, equivalently, \eqref{lim0} holds.
\end{enumerate}
\end{proposition}

\begin{proof}
Since $A$ is normal, it follows that
%To see the equivalence of (i)--(iv) observe that if $A$ is normal
%then
$D_{A^{*n}}=D_{A^{n}}$ for all $n\in\dN_0$. Hence the identities
\begin{equation}\label{cont}
\left(\sH^c_\Sigma\right)^\perp=\bigcap_{n=0}^\infty
\ker(D_{A^*}A^{*n})=\bigcap_{n=1}^\infty \ker D_{A^{*n}},
\end{equation}
\begin{equation}\label{obs}
\left(\sH^o_\Sigma\right)^\perp=\bigcap_{n=0}^\infty
\ker(D_{A}A^{n})=\bigcap_{n=1}^\infty \ker D_{A^{n}}
\end{equation}
imply that
$\left(\sH^c_\Sigma\right)^\perp=\left(\sH^o_\Sigma\right)^\perp$,
or equivalently, $\sH^c_\Sigma=\sH^o_\Sigma$. This identity implies
the equivalence of (i), (ii), (iii), and (iv).
%The equivalence of
%(iv) and (v) is a general fact for not necessarily normal
%contractions.
\end{proof}

The following corollary is based on the fact that a contraction $A$
is completely nonunitary if and only if the corresponding system
$\Sigma$ in \eqref{SIGMA}, \eqref{blockmatrix1}  is minimal.

\begin{corollary}\label{stablecor}
Let $A$ be a normal contraction. Then the statements (i)--(iii) in
Proposition \ref{nstable} and the statements (i)--(iv) in
Proposition \ref{nnstable} are all equivalent.
\end{corollary}

\subsection{Parametrization of block operators}

For a proof and some history of the following theorem, see
\cite{AHS2}.

\begin{theorem} \label{ParContr}
Let $\sM, \sN$, $\sH$, and $\sK$ be Hilbert spaces. The operator
matrix $T$ in \eqref{abcd0}  is a contraction if and only if $T$ is
of the form
\begin{equation}\label{twee}
T=
\begin{pmatrix} -KA^*M+D_{K^*}XD_{M} &KD_{A} \cr
D_{A^*}M&A\end{pmatrix},
\end{equation}
where $A \in \bL(\sH,\sK)$, $M\in\bL(\sM,\sD_{A^*})$,
$K\in\bL(\sD_{A},\sN)$, and $X\in\bL(\sD_{M},\sD_{K^*})$ are
contractions, all uniquely determined by $T$. Furthermore, the
following equality holds for all $h \in \sM$,$f \in \sH$:
\begin{equation} \label{contr}
\begin{split}
\left\|\begin{pmatrix}h\cr f\end{pmatrix}\right\|^2
&-\left\|\begin{pmatrix} -KA^*M+D_{K^*}XD_{M} &KD_{A} \cr
D_{A^*}M&A\end{pmatrix}
\begin{pmatrix}h\cr f\end{pmatrix}  \right\|^2 \\
&=\|D_K(D_{A}f-A^*Mh)-K^*XD_{M}h\|^2+\|D_{X}D_{M}h\|^2.
\end{split}
\end{equation}
\end{theorem}
\begin{corollary}\label{Corvaasa}
\label{iso} Let $A\in\bL(\sH,\sK)$ be a contraction. Assume that
$K\in\bL(\sD_A,\sN)$, $M\in\bL(\sM,\sD_{A^*})$, and
$X\in\bL(\sD_M,\sD_{K^*})$ are contractions. Then the operator $T$
in \eqref{twee} is:
\begin{enumerate}
\def\labelenumi{\rm (\roman{enumi})}
\item isometric if and only if $D_XD_M=0$ and $D_KD_A=0$;
\item co-isometric if and only if $D_{X^*}D_{K^*}=0$ and
$D_{M^*}D_{A^*}=0$.
\end{enumerate}
\end{corollary}

Let $\tau=\{T;\sM,\sN,\sH\}$ be a passive system and let
\eqref{twee} be the representation of the block operator $T$ in
\eqref{abcd0}. Define for $\lambda\in\dD$ the following
operator-valued holomorphic functions
 \begin{equation}
 \label{fi}
\f(\lambda):=\begin{pmatrix} -D_XD_{M} \cr
D_{K}\Phi_{A^*}(\lambda)M-K^*XD_{M}
\end{pmatrix}:\sM\to \begin{pmatrix}\sD_M\\
\sD_K\end{pmatrix},
\end{equation}
and
\begin{equation}
\label{psi}
 \psi(\lambda):=\begin{pmatrix}  D_{K^*}D_{X^*}&K\Phi_{A^*}(\lambda)
D_{M^*}-D_{K^*}XM^* \end{pmatrix}:
\begin{pmatrix}\sD_{K^*}\\
\sD_{M^*}\end{pmatrix} \to\sN.
\end{equation}

\begin{theorem}[\cite{AHS2}]\label{Schur}
Let $\tau=\{T;\sM,\sN,\sH\}$ be a passive system and let
\eqref{twee} be the representation of the block operator $T$ in
\eqref{abcd0}. Then the transfer function $\Theta (\lambda)$ of
$\tau$ and the characteristic function $\Phi_{A^*}(\lambda)$ of
$A^*$ (see \eqref{CHARFUNC1}) are connected via
\begin{equation}
\label{exptran}
 \Theta (\lambda)=K\Phi_{A^*}(\lambda)M+D_{K^*}XD_{M},\quad
\lambda\in\dD;
\end{equation}
in particular, $\Theta (\lambda) \in {\bf S}(\sM,\sN)$. In addition,
the identities
\begin{equation}
\label{deffunc1} \left\|D_{\Theta (\lambda)}h\right\|^2 =
\left\|D_{\Phi_{A^*}(\lambda)}Mh\right\|^2  +\|\varphi(\lambda)
h\|^2,
  \quad  h\in\sM,
\end{equation}
\begin{equation}
\label{deffunc2} \left\|D_{\Theta^* (\lambda)}g\right\|^2 =
\left\|D_{\Phi_{A}(\ovl\lambda)}K^*g\right\|^2 + \|\psi^*(\lambda)
g\|^2, \quad g\in\sN,
\end{equation}
hold and the functions $\varphi(\lambda)$ and $\psi(\lambda)$ in
\eqref{fi} and \eqref{psi} are Schur functions.
\end{theorem}

\section{Passive systems with a normal main operator}

Let $\tau$ be a passive system  of the form \eqref{abcd0}. If its
main operator $A$ is normal, then  many properties of $\tau$ and its
transfer function   simplify.

\subsection{Basic properties}

The controllable and observable subspaces of the passive system in
\eqref{abcd0} are defined in \eqref{CO}. Let the block matrix $T$
have the parametrization   \eqref{twee}, so that  $A^nB=A^nD_{A^*}M$
and $A^{*n}C^*=A^{*n} D_AK^*$.  If, in addition,  $A$ is normal it
follows that $D_{A^*}=D_A$ and then \eqref{comm+} implies
\[
 A^nB=D_A A^n M, \quad A^{*n}C^*=A^{*n} D_AK^*.
\]
Hence, if $A$ is normal, then $\sH^c$ and $\sH^o$ have the form:
\begin{equation}
\label{CON1}
 \sH^c=\cspan\{\,D_{A} A^{n} M \sM:\,n \in \dN_0 \},
\quad
 \sH^o=\cspan\{\,D_{A} A^{*n} K^* \sN:\, n \in \dN_0 \},
\end{equation}
or, equivalently,
\begin{equation}\label{CON3}
 (\sH^c)^\perp=\bigcap\limits_{n=0}^\infty\ker(M^*A^{*n}D_{A}),
  \qquad
 (\sH^o)^\perp=\bigcap\limits_{n=0}^\infty\ker(KA^nD_{A}),
\end{equation}
%These representations  give rise to the following definitions. Let
Let the subspaces $\sH^c_N$ and $\sH^o_N$ be defined by
\begin{equation}\label{SPn}
\sH^c_N=\cspan\left\{\,A^nM\sM:\, n \in \dN_0 \, \right\}, \quad
\sH^o_N=\cspan\left\{\,A^{*n}K^*\sN:\, n \in \dN_0 \, \right\},
\end{equation}
or, equivalently, by
\begin{equation}\label{CON4}
 (\sH_N^c)^\perp=\bigcap\limits_{n=0}^\infty\ker(M^*A^{*n}),
 \qquad
 (\sH_N^o)^\perp=\bigcap\limits_{n=0}^\infty\ker(KA^n).
\end{equation}

\begin{lemma}\label{normal0}
Let $\tau=\{T;\sM,\sN,\sH\}$ be a passive system with $T$   of the
form \eqref{twee} with some contractions $A \in \bL(\sH,\sK)$,
$M\in\bL(\sM,\sD_{A^*})$, $K\in\bL(\sD_{A},\sN)$, and
$X\in\bL(\sD_{M},\sD_{K^*})$. Assume that $A$ is normal.
\begin{enumerate}
\def\labelenumi{\rm (\roman{enumi})}
\item If $\sH^c_N$ is invariant under $A^*$, then
$\sH\ominus\sH_N^c\subset \sH\ominus\sH^c$.
\item If $\sH_N^o$ be invariant under $A$, then
$\sH\ominus\sH_N^o\subset \sH\ominus\sH^o$.
\end{enumerate}
\end{lemma}

\begin{proof}
(i) Assume that $\sH^c_N$ is invariant under $A^*$ or, equivalently,
that $\sH\ominus\sH^c_N$ is invariant under $A$. Hence, if
$f\in\sH\ominus\sH^c_N$ then $f$ and $Af$ both belong to
$\ker(M^*A^{*n})$ for all $n\in \dN_0$. Thus, in particular,
$D^2_Af=(I-A^*A)f\in\ker M^*$. Moreover, if $p(t)$ is a polynomial
then $p(D^2_A)f\in\ker M^*$. Since there exists a sequence of
polynomials $\{p_m(t)\}_{m=1}^\infty$ such that the sequence
$\{p_m(D^2_A)\}$ converges uniformly  to $D_A$, it follows that
$D_Af\in\ker M^*$. Furthermore, the sequence $\{p_m(D^2_A)A^{*n}\}$
converges uniformly to $D_AA^{*n}$ for all $n \in \mathbb{N}$. Since
$p_m(D^2_A)A^{*n}f\in \ker M^*$ for all $n\in\dN_0$, one concludes
that $D_AA^{*n} f\in\ker M^*$ for all $n \in \mathbb{N}_0$. It
follows that
 \[
 \sH\ominus\sH_N^c\subset \bigcap\limits_{n=0}^\infty\ker(M^* D_A A^{*n})
 =\sH\ominus\cspan\left\{\,D_AA^nM\sN:\,
n \in \dN_0 \,\right\}=\sH\ominus\sH^c.
\]

(ii) The proof of (ii) is similar to the proof of (i).
\end{proof}

\begin{proposition}\label{normal1}
Let $\tau=\{T;\sM,\sN,\sH\}$ be a passive system where $T$ is of the
form \eqref{twee} with contractions $A \in \bL(\sH,\sK)$,
$M\in\bL(\sM,\sD_{A^*})$, $K\in\bL(\sD_{A},\sN)$, and
$X\in\bL(\sD_{M},\sD_{K^*})$. Assume that $A$ is normal.
\begin{enumerate}
\def\labelenumi{\rm (\roman{enumi})}
\item $\tau$ is controllable if and only if
\begin{equation}\label{COa}
 \ker D_A=\{0\} \quad \text{and}\quad
 \ran D_A \cap \left(\sH\ominus\sH^c_N\right)=\{0\}.
\end{equation}
In particular, if $\ker D_A=\{0\}$ and $\sH^c_N=\sH$, then   $\tau$
is controllable;   if $\tau$ is controllable and $\sH_N^c$ is
invariant under $A^*$, then $\ker D_A=\{0\}$ and $\sH^c_N=\sH$.

\item $\tau$ is observable if and only if
\begin{equation}\label{COb}
 \ker D_A=\{0\} \quad \text{and}\quad
 \ran D_A \cap \left(\sH\ominus\sH^o_N\right)=\{0\};
\end{equation}
In particular, if $\ker D_A=\{0\}$ and $\sH^o_N=\sH$, then   $\tau$
is observable;   if $\tau$ is observable and $\sH_N^o$ is invariant
under $A $, then $\ker D_A=\{0\}$ and $\sH^o_N=\sH$.

\item $\tau$ is simple if and only if
\begin{equation}\label{COc}
 \ker D_A=\{0\} \quad \text{and}\quad
 \ran D_A \cap \sH\ominus\left(\sH^c_N+\sH^o_N\right)=\{0\}.
\end{equation}
In particular, if $\ker D_A=\{0\}$ and $\sH=\clos
\{\sH_N^c+\sH_N^o\}$, then $\tau$ is simple; if $\tau$ is simple,
$\sH_N^c$ is invariant under $A^*$, and $\sH_N^o$ is invariant under
$A$, then $\ker D_A=\{0\}$ and $\sH=\clos \{\sH_N^c+\sH_N^o\}$.

\item $\tau$ is minimal if and only if
\begin{equation}\label{COa1}
\begin{split}
 \ker D_A&=\{0\}, \quad
 \ran D_A \cap \left(\sH\ominus\sH^c_N\right)=\{0\}, \\
 &\hspace{3.5cm}\text{and}\quad  \ran D_A \cap
 \left(\sH\ominus\sH^o_N\right)=\{0\}.
\end{split}
\end{equation}
%the conditions in (i) and (ii) hold.
In particular, if $\ker D_A=\{0\}$ and $\sH=\sH^c_N=\sH^o_N$, then
$\tau$ is minimal; if $\tau$ is minimal, $\sH_N^c$ is invariant
under $A^*$, and $\sH_N^o$ is invariant under $A$, then $\ker
D_A=\{0\}$ and $\sH=\sH^c_N=\sH^o_N$.
\end{enumerate}
\end{proposition}

\begin{proof}
(i) Assume that \eqref{COa} holds. Let $f\in(\sH^c)^\perp$. It
follows from \eqref{CON3} and \eqref{CON4} that $D_Af \in
\bigcap_{n=0}^\infty\ker(M^*A^{*n})=(\sH_N^c)^\perp$. The second
condition in \eqref{COa} shows that $D_Af=0$ and the first condition
in \eqref{COa} yields $f=0$. Therefore, $\sH^c=\sH$ and $\tau$ is
controllable.

Now assume that $\tau$ is controllable, i.e. $\sH^c=\sH$.  Then
\eqref{CON3} implies that $\ker D_A=\{0\}$. Furthermore, if $D_Af
\in (\sH_N^c)^\perp$, then   \eqref{CON4} implies that $D_Af \in
\bigcap_{n=0}^\infty\ker(M^*A^{*n})$. By  \eqref{CON3} this leads to
$f\in (\sH^c)^\perp$ and hence $f=0$ by controllability of $\tau$.
This shows that \eqref{COa} is satisfied.

If $\ker D_A=\{0\}$ and $\sH^c_N=\sH$, then  \eqref{COa} is
satisfied. It follows that $\tau$ is controllable.

If $\sH^c_N$ is invariant under $A^*$, then
$\sH\ominus\sH_N^c\subset \sH\ominus\sH^c$ by Lemma \ref{normal0}.
Hence, if in addition, $\tau$ is controllable, it follows that
$\sH^c_N=\sH$; moreover, it follows that $\ker D_A=\{0\}$.

(ii) The proof is completely analogous to the proof for part (i).

(iii) If $\tau$ is simple then it immediately follows from
\eqref{CON1} that $\ker D_A=\{0\}$. Moreover, it is
clear from \eqref{CON3} and \eqref{CON4} that $D_Af\in
\sH\ominus\left(\sH^c_N+\sH^o_N\right)$ if and only if $f\in
\sH\ominus\left(\sH^c+\sH^o\right)$. Now the statement is obtained
as in part (i).

(iv) This is obvious from the definition of minimality.
\end{proof}

\begin{corollary}\label{Norcoro}
Let the main operator $A$ of the passive system
$\tau=\{T;\sM,\sN,\sH\}$ be normal and let the system be simple.
Then the system $\tau$ is strongly stable and strongly co-stable.
\end{corollary}

\begin{proof}
Since $\tau$
%=\{T;\sM,\sN,\sH\}$
is simple and $A$ is normal, Proposition \ref{normal1} shows that
$\ker D_A=\{0\}$ or, equivalently, that the contraction $A$ is
strict. Hence Lemma~\ref{nstable} implies that $A\in C_{00}$.
Therefore $\tau$ is strongly stable and strongly co-stable.
\end{proof}

\subsection{Defect functions}

Associated with $\Theta(\lambda) \in \mathbf{S}(\sM,\sN)$ are the
right and left \textit{defect functions} (or \textit{spectral
factors}) $\f_\Theta(\lambda)$ and $\psi_\Theta(\lambda)$, which
satisfy
\begin{equation}\label{einz}
 \varphi^*_\Theta(\xi)\varphi_\Theta(\xi)\le I_\sM
-\Theta^*(\xi)\Theta(\xi),
\quad
 \psi_\Theta(\xi)\psi^*_\Theta(\xi)\le I_\sN-
\Theta(\xi)\Theta^*(\xi),
\end{equation}
almost everywhere on $\dT$. These operator-valued Schur functions
are (up to a constant unitary factor) uniquely determined by the
following maximality property: if $\wt\varphi(\lambda)$ and
$\wt\psi(\lambda)$ are operator-valued Schur functions for which
\begin{equation}\label{zwei}
\wt\varphi^*(\xi)\wt\varphi(\xi)\le I_\sM -\Theta^*(\xi)\Theta(\xi),
\quad \wt\psi(\xi)\wt\psi^*(\xi)\le I_\sM -\Theta(\xi)\Theta^*(\xi),
\end{equation}
then they are dominated by $\f_\Theta(\lambda)$ and
$\psi_\Theta(\lambda)$  in the following sense:
\begin{equation}\label{drei}
\wt\varphi^*(\xi)\wt\varphi(\xi)\le
\varphi^*_\Theta(\xi)\varphi_\Theta(\xi), \quad
\wt\psi(\xi)\wt\psi^*(\xi)\le \psi_\Theta(\xi)\psi^*_\Theta(\xi),
\end{equation}
almost everywhere on the unit circle $\dT$; cf. \cite{BC},
\cite{BD}, \cite{BDFK1}, \cite{BDFK2}, \cite{DR}.

Note that it follows from Theorem \ref{Schur} that the functions
$\varphi(\lambda)$ and $\psi(\lambda)$ satisfy the inequalities
\begin{equation}\label{uksi}
\varphi(\xi)^*\varphi(\xi) \le
\varphi^*_\Theta(\xi)\varphi_\Theta(\xi),  \quad
\psi(\xi)\psi^*(\xi) \le \psi_\Theta(\xi)\psi^*_\Theta(\xi),
\end{equation}
for almost all $\xi \in \dT$.

\begin{proposition}\label{NorIN}
Let $\tau=\{T;\sM,\sN,\sH\}$ be a passive system with a normal main
operator $A$ and let $\Theta (\lambda)\in{\bf S}(\sM,\sN)$ be its
transfer function. If $\tau$ is simple and
$\varphi_\Theta(\lambda)=0$ ($\psi_\Theta(\lambda)=0$), then
$\Theta(\lambda)$ is inner (co-inner, respectively).
\end{proposition}

\begin{proof} By Corollary~\ref{Norcoro} one has $A\in C_{00}$ and, in
particular, $A$ is completely non-unitary. Therefore,
$\Phi_A(\lambda)$ and $\Phi_{A^*}(\lambda)$ are bi-inner. On the
other hand, if $\varphi_\Theta(\lambda)=0$
($\psi_\Theta(\lambda)=0$), then \eqref{uksi} shows that
$\varphi(\xi)=0$ ($\psi(\xi)=0$) for almost all $\xi\in\dT$. Now
\eqref{deffunc1} (\eqref{deffunc2}, respectively) yields that
$D_{\Theta(\xi)}=0$ ($D_{\Theta^*(\xi)}=0$) almost everywhere on
$\dT$, i.e. $\Theta(\lambda)$ is inner (co-inner).
\end{proof}

\subsection{Unitary similarity}

Recall that two passive systems
$\tau_j=\left\{T_j;\sN,\sN,\sH_j\right\}$, $j=1,2$, are said to be
unitarily similar if there is a unitary operator $U:\sH_1\to \sH_2$,
such that \eqref{unisim} holds.
%\begin{equation}\label{unisim}
% A_2=UA_1U^*,\quad  B_2=UB_1, \quad  C_2=C_1U^*, \mbox{ and } D_1=D_2.
%\end{equation}
In particular, in this case the spectra of the corresponding main
operators $A_1$ and $A_2$ coincide. It is clear that if the systems
$\tau_1$ and $\tau_2$ are unitarily similar then they have the same
transfer function. However, two minimal passive systems $\tau_1$ and
$\tau_2$ with the same transfer function $\Theta(\lambda)$ are in
general not unitarily similar; such systems are only weakly similar
as shown in D.Z.~Arov~\cite{A}, see \eqref{pseudosim}. In the case
of passive systems with  normal main operators the following
sufficient spectral-theoretic condition can be established.

\begin{theorem}\label{NUS}
Let $\tau_1=\left\{T_1;\sM,\sN, \sH_{1}\right\}$ and
$\tau_2=\left\{T_2;\sM,\sN, \sH_{2}\right\}$ be two minimal passive
systems whose transfer functions coincide in some neighborhood of
zero. Let the main operator $A_k$ be normal and let $C_k=SB_k^*$,
$k=1,2$, with $S$ bounded and injective. Then, if the spectrum
$\sigma(A_k)$ of $A_k$, $k=1,2$, does not contain interior points
and $\rho(A_1)\cap \rho(A_2)$ is a connected set in $\dC$, the
systems $\tau_1$ and $\tau_2$ are unitarily similar.
\end{theorem}

\begin{proof}
Assume that the transfer functions $\Theta_{1}(\lambda)$ and
$\Theta_{2}(\lambda)$ of $\tau_1$ and $\tau_2$  coincide in some
neighborhood of zero. Since $\Theta_{ 1}(\lambda)$ and $\Theta_{
2}(\lambda)$ are holomorphic on $\dD$ it follows that $\Theta_{
1}(\lambda)=\Theta_{ 2}(\lambda)$ for all $\lambda \in \dD$. The
definition \eqref{TrFu} implies that $D_1=\Theta_{ 1}(0)=\Theta_{
2}(0)=D_2$ and that
\[
\sum_{m=0}^\infty \lambda^mC_1A^m_{1}B_1= \sum_{m=0}^\infty
\lambda^mC_2A^m_{2}B_2,\quad \lambda \in \dD.
\]
Since $C_k=SB^*_k$, $k=1,2$, where $S$ is bounded and injective, the
previous equality yields
\begin{equation}\label{degn0}
 B^*_1 A^m_ {1}B_1=B^*_2 A^m_{2} B_2,\quad m\in \mathbb{N}_0.
\end{equation}
%Now define
%\[
% Z_0\left(\sum\limits_{j=0}^mA^j_{1}B_1u_j\right)=
% \sum\limits_{j=0}^mA^j_{2}B_2u_j, \quad u_0, u_1,\ldots,u_m \in \sM,
% \quad m \in \dN_0.
%\]
%Then $Z_0$ is linear, $\dom Z_0=\lin\{\,A_{1}^{n}B_1\sM:\,n \in
%\dN_0 \}$, $\ran Z_0=\lin\{\,A_{2}^{n}B_2\sM:\,n \in \dN_0 \}$, and
%in particular, $Z_0$ is densely defined ($\sH_1^c=\sH_1$). Hence,
%$\mul Z_0^*=(\dom Z_0^{**})^\perp=\{0\}$ and thus the adjoint
%$Z_0^*$ is an operator. It follows from \eqref{degn0} that $Z_0^*$
%satisfies
%\begin{equation}\label{Z0*}
% Z_0^*\left(\sum\limits_{j=0}^m A^{*j}_{2}B_2 v_j\right)=
% \sum\limits_{j=0}^m A^{*j}_{1}B_1 v_j, \quad v_0, v_1,\ldots,v_m \in \sM,
% \quad m \in \dN_0.
%\end{equation}
%In particular, by observability assumptions (note that
%$C_k^*=B_kS^*$), $Z_0^*$ is densely defined with dense range.
%Therefore, $Z_0^{**}$ is an operator, i.e., $Z_0$ is closable
%(compare \cite{Arov}).
Now define the relation $Z_0$ by
\begin{equation}\label{Z0}
 Z_0=\left\{\, \left\{\sum\limits_{j=0}^mA^j_{1}B_1u_j,
 \sum\limits_{j=0}^mA^j_{2}B_2u_j \right\} :\, u_0, u_1,\ldots,u_m \in \sM, \,\,m \in
 \dN_0\, \right\}.
\end{equation}
Clearly $Z_0$ is linear and
\[
\dom Z_0=\lin\{\,A_{1}^{n}B_1\sM:\,n \in \dN_0 \}, \quad \ran
Z_0=\lin\{\,A_{2}^{n}B_2\sM:\,n \in \dN_0 \}.
\]
Furthermore, it follows from \eqref{degn0} that
\begin{equation}\label{Z0*}
 \left\{ \sum\limits_{j=0}^m A^{*j}_{2}B_2 v_j, \sum\limits_{j=0}^m A^{*j}_{1}B_1
 v_j \right\} \in Z_0^*, \quad v_1,\ldots,v_m \in \sM,
 \quad m \in \dN_0,
\end{equation}
so that
\[
\lin\{\,A_{2}^{*n}B_2\sM:\,n \in \dN_0 \} \subset \dom Z_0^*, \quad
\lin\{\,A_{1}^{*n}B_1\sM:\,n \in \dN_0 \} \subset \ran Z_0^*.
\]
Due to the controllability and observability conditions (note that
$C_k^*=B_kS^*$), it follows from \eqref{Z0} and \eqref{Z0*} that
both $Z_0$ and $Z_0^*$ have dense domains and dense ranges. In
particular, $Z_0$ and $Z_0^*$ are (graphs of) operators, and, in
fact, $\mul Z_0^{**}=(\dom Z_0^*)^\perp$ implies that $Z_0$ is a
closable operator, i.e., its closure $Z_0^{**}$ is (the graph of) an
operator; cf. \cite{Arov}.

Next it is shown that under the assumptions on the main operators
$A_k$, $k=1,2$, the mapping $Z_0$ becomes isometric. Since $A$ is
contractive the spectrum $\sigma(A_k)$ is a compact subset of the
closed unit disk. The  union $\sigma(A_1)\cup \sigma(A_2)$ is also
compact and, in addition, does not have interior points. Indeed,
this follows immediately from the fact that the sets $\sigma(A_1)$
and $\sigma(A_2)$ are closed and do not have interior points.
Furthermore, by assumption $\dC\setminus(\sigma(A_1)\cup
\sigma(A_2))=\rho(A_1)\cap \rho(A_2)$ is connected. Therefore,
according to Mergelyan's theorem (see e.g.
\cite[Theorem~20.5]{Rudin}) every continuous complex-valued function
on $\sigma(A_1)\cup \sigma(A_2)$ can be uniformly approximated on
$\sigma(A_1)\cup \sigma(A_2)$ by complex polynomials. Since for
every $n,m\in \dN_0$ the function $f_{n,m}(z)=\ovl{z}^nz^m$ is
continuous on $\dC$, there exists a sequence
$\{P_j^{n,m}(z):\,j\in\dN_0\}$ of polynomials converging uniformly
on $\sigma(A_1)\cup \sigma(A_2)$ to $f_{n,m}(z)$. It follows from
\eqref{degn0} that for every $n,k,j\in\dN_0$ one has
\begin{equation}\label{degn1}
 B^*_1 P_j^{n,m}(A_1) B_1=B^*_2 P_j^{n,m}(A_2) B_2.
\end{equation}
The functional calculus for normal operators shows that
$f_{n,m}(A_k)=A_k^{*n}A_k^m$, $k=1,2$, and therefore  taking strong
limits in \eqref{degn1} yields
\begin{equation}\label{degn2}
 B^*_1 A_1^{*n}A_1^m B_1=B^*_2 A_2^{*n}A_2^m B_2,
 \quad m,n \in \dN_0.
\end{equation}
These identities imply that
\[
 \left\|\sum_{j=0}^mA^j_{1}B_1u_j \right\|^2 =
\left\|\sum_{j=0}^mA^j_{2}B_2u_j \right\|^2, \quad u_0,
u_1,\ldots,u_m \in \sM, \quad m \in \dN_0,
\]
and, therefore, the operator $Z_0$ in \eqref{Z0} is isometric. Since
$Z_0$ is densely defined with dense range, its closure $Z$ is
unitary. The identities $ZA_1=A_2Z$ and $ZB_1=B_2$ are immediate
from \eqref{Z0}, while \eqref{Z0*} shows that $Z_0^*B_2=B_1$ which
gives the identity $C_2Z=C_1$. Therefore, the systems $\tau_1$ and
$\tau_2$ are unitarily similar; cf. \eqref{unisim}.
\end{proof}

%In particular, the conditions in Theorem~\ref{NUS} imply that the
%spectra of $A_1$ and $A_2$ must coincide. The following two
%corollaries are immediate from this theorem.

\begin{corollary}\label{US1}
Let $\tau_1=\left\{T_1;\sN,\sN, \sH_{1}\right\}$ and
$\tau_2=\left\{T_2;\sN,\sN, \sH_{2}\right\}$ be two minimal passive
systems such that $A_k$ is selfadjoint ($A_k=A_k^*$) or
skew-symmetric ($A_k=-A_k^*$) and $C_k=SB_k^*$, $k=1,2$, with $S$
bounded and injective. Then $\tau_1$ and $\tau_2$ are unitarily
similar if and only if their transfer functions coincide in some
neighborhood of zero.
\end{corollary}

\begin{corollary}\label{US12}
Let $\tau_1=\left\{T_1;\sN,\sN, \sH_{1}\right\}$ and
$\tau_2=\left\{T_2;\sN,\sN, \sH_{2}\right\}$ be two minimal passive
systems such that $A_k$ is normal and has a discrete spectrum,
and $C_k=SB_k^*$, $k=1,2$, with $S$ bounded and injective.
Then $\tau_1$ and $\tau_2$ are unitarily
similar if and only if their transfer functions coincide in some
neighborhood of zero.
\end{corollary}

\begin{corollary}\label{US2}
Let $\tau_1=\left\{T_1;\sN,\sN, \sH_{1}\right\}$ and
$\tau_2=\left\{T_2;\sN,\sN, \sH_{2}\right\}$ be two minimal passive
systems with a finite-dimensional state space $\sH_k$ such that
$A_k$ is normal and $C_k=SB_k^*$, $k=1,2$, with $S$ bounded and
injective. Then $\tau_1$ and $\tau_2$ are unitarily similar if and
only if their transfer functions coincide in some neighborhood of
zero.
\end{corollary}

\begin{remark}
(i) The proof of Theorem~\ref{NUS} uses that fact that the operators
$f_{n,m}(A)=A^{*n}A^m$, $m,n\in\dN_0$, can be approximated by a
sequence of polynomials in $A$. If, in particular, the adjoint $A^*$
of a bounded operator $A$ can be approximated by a sequence
$P_n(A)$, $n\in\dN_0$, of polynomials in $A$, i.e.,
\begin{equation}
\label{p0}
 A^*=s-\lim\limits_{n\to\infty} P_n(A),
\end{equation}
then the same is true for all of the operators
$f_{n,m}(A)=A^{*n}A^m$, $m,n\in\dN_0$. By taking strong limits in
$AP_n(A)=P_n(A)A$ one obtains from \eqref{p0} the identity
$AA^*=A^*A$. Therefore, the condition \eqref{p0} implies that $A$ is
a normal operator.

(ii) If $A$ is a normal operator, then $A^*=f(A)$ with
$f(z)=\ovl{z}$ by the functional calculus for normal operators. The
function $f(z)$ does not satisfy the Cauchy-Riemann equations, so it
is nowhere holomorphic. Consequently, if $\sigma(A)$ has interior
points, the adjoint $A^*$ cannot satisfy the condition \eqref{p0},
as one would get a uniform approximation for $f(z)$ on $\sigma(A)$
via polynomials $P_n(z)$.

(iii) If $A$ is a normal operator on a finite-dimensional space, then
it has $n=\dim \sH$ eigenvalues and it is unitarily similar to a
diagonal matrix. Therefore, if $A$ has $d$ nonreal eigenvalues then
by standard interpolation one finds a polynomial $Q$
(say, of degree at most $d-1$ when using only the nonreal spectral points)
such that $A^*=Q(A)$.
%; see e.g. \cite{HJ}.
If there are two normal operators $A_1$ and $A_2$ on
$\sH_k$, $n_k=\dim \sH_k <\infty$, $k=1,2$, then together they have at
most $n_1+n_2$ different nonreal eigenvalues and one can find a
polynomial $P$ (of degree at most $n_1+n_2-1$) such that
$A_1^*=P(A_1)$ and $A_2^*=P(A_2)$. Then
$f_{n,m}(A_k)=A_k^{*n}A_k^m=P(A_k)^nA_k^m$, $n,m\in\dN_0$, is also a
polynomial in $A_k$, $k=1,2$. So, in the proof of Theorem~\ref{NUS}
no limit procedure is needed in the case of finite-dimensional state
spaces.

(iv) Finally, note that the criterion for unitary similarity
of minimal passive systems with the same transfer function which has
been established in \cite{ArNu1} is essentially of different nature
than the above spectral theoretical sufficient condition in
Theorem~\ref{NUS}.
\end{remark}

\section{Passive quasi-selfadjoint systems}

\subsection{Quasi-selfadjoint contractions and associated passive systems}

Let $\cH$ be a Hilbert space. A linear operator $T\in\bL(\cH)$ is
said to be a \textit{quasi-selfadjoint contraction} ($qsc$-operator
for short) if
\[
\dom T=\cH,\; \|T\|\le 1, \mbox{ and } \ker(T-T^*)\ne \{0\}.
\]
The next theorem is a consequence of Theorem \ref{ParContr}; see
\cite{AHS1}.

\begin{theorem}\label{qscex}
Let $T$ be a $qsc$-operator in the Hilbert space $\cH$ and let $\sN$
be a subspace in $\cH$ such that $\ran(T-T^*)\subset \sN$. Then with
respect to the decomposition $\cH=\sN\oplus\sH$, where
$\sH=\cH\ominus\sN$, the operator $T$ has the following block form
\begin{equation}\label{MCE}
T=\begin{pmatrix} -KAK^*+D_{K^*}XD_{K^*} &KD_{A}\cr D_{A}K^* &
A\end{pmatrix}:\begin{pmatrix}\sN\\ \sH\end{pmatrix} \to
\begin{pmatrix}\sN\\ \sH\end{pmatrix},
\end{equation}
where $A=P_\sH T\uphar\sH$ is a selfadjoint contraction and
$K\in\bL(\sD_A,\sN)$, $X\in \bL(\sD_{K^*})$ are contractions.
\end{theorem}

The  system $\tau=\left\{T;\sN,\sN,\sH\right\}$ is said to be
\textit{passive quasi-selfadjoint} ($\tau$ is a $pqs$-system for
short) if  $T$ in \eqref{abcd} is a contraction  and if ${\rm
ran\,}(T-T^*)\subset\mathfrak N$.  It follows  that   $T$ is a
$qsc$-operator in $\sN\oplus\sH$ and that $A=A^*$ and $C=B^*$.
Moreover, according to Theorem \ref{qscex},  $B$, $C,$ and $D$ have
the form
\begin{equation}
\label{bcd}
 B=D_{A}K^*,\quad C=KD_A,\quad
D=-KAK^*+D_{K^*}XD_{K^*},
\end{equation}
where $K\in\bL(\sD_{A},\sN)$ and $X\in\bL(\sD_K^*)$ are
contractions. For a $pqs$-system the controllable and observable
subspaces  coincide, see \eqref{CON1}:
\begin{equation}\label{neww}
 \sH^c = \sH^o  =\cspan\left\{\,D_AA^nK^*\sN:\, n \in \dN_0
\,\right\}\subset\sD_A.
\end{equation}

\subsection{Minimal representations of $pqs$-systems and unitary similarity}

A $pqs$-system can always be reduced to a minimal $pqs$-system.

\begin{proposition}\label{simple1}
Let $\tau=\left\{T;\sN,\sN,\sH\right\}$ be a $pqs$-system of the
form \eqref{abcd} and let $B,C,$ and $D$ be given by \eqref{bcd}
with some contractions $K$ and $X$. Define the system
\begin{equation}\label{minsys1}
 \tau_s=\{T_s;\sN,\sN,\sH^s\},
\end{equation}
where the subspace $\sH^s$ is given by
\begin{equation}\label{SP}
\sH^s=\cspan\left\{\,A^nK^*\sN:\, n \in \dN_0 \, \right\},
\end{equation}
and where the operator $T_s$ is given by
\begin{equation}\label{abcds}
T_s=\begin{pmatrix} D&C\uphar\sH^s\cr B&A\uphar\sH^s\end{pmatrix}:
\begin{pmatrix}\sN\\ \sH^s\end{pmatrix}\to
\begin{pmatrix}\sN \\ \sH^s\end{pmatrix}.
\end{equation}
Then $\tau_s$ is a minimal $pqs$-system and the transfer functions
of the systems $\tau$ and $\tau^s$ coincide. Moreover, the system
$\tau$ is minimal if and only if
\begin{enumerate}\def\labelenumi{\rm (\roman{enumi})}
\item
$\left\|Af\right\|<\left\|f\right\|$ for all $f\in
\sH\backslash\{0\},$
\item $\sH^s=\sH$.
\end{enumerate}
In this case the system $\tau$ is strongly stable and strongly
co-stable.
\end{proposition}

\begin{proof}

The subspace $\sH^s$ in \eqref{SP} reduces $A$ and therefore it also
reduces $D_A=(I_\sH-A^2)^{1/2}$. Furthermore, $\cran
K^*\subset\sH^s$. Let $ A_s=A\uphar\sH^s$, then $D_{
A_s}=D_A\uphar\sH^s $ and, hence, $D_AK^*=D_{A_s}K^*$. Define the
operator $C_s$ by
\[
 C_s=C\uphar\sH^s=KD_{A_s}.
\]
Then $T_s$ in \eqref{abcds} is a $qsc$-operator in $\sN\oplus\sH^s$.
Since $\ran K^*\subset \sD_A\cap \sH^s$, one has $\sD_{A_s}=\sH^s$.
Now the construction shows that the system $\tau^s$ in
\eqref{minsys1} is minimal. Clearly, the transfer functions of
$\tau$ and $\tau^s$ coincide.

As to the minimality of $\tau$ observe that $\sH^s=\sH^c_N=\sH^o_N$,
since $A=A^*$; see \eqref{SPn}. Hence, the
characteristic properties (i) and (ii) for minimality of a
$pqs$-system $\tau$ are obtained from Proposition~\ref{normal1}.

The last statement holds by Corollary~\ref{Norcoro}.
\end{proof}

%\begin{corollary}\label{coro}
%\label{col1} A minimal $pqs$-system is strongly stable and strongly
%co-stable.
%\end{corollary}
%\begin{proof}
%The statement holds by Corollary~\ref{Norcoro}.
%\end{proof}

It is a consequence of Theorem~\ref{NUS} that within the class of
$pqs$-systems the following unitary similarity criterion holds; see
Corollary~\ref{US1}.

\begin{proposition}\label{US}
Let $\tau_1=\left\{T_1;\sN,\sN, \sH_{1}\right\}$ and
$\tau_2=\left\{T_2;\sN,\sN, \sH_{2}\right\}$ be two minimal
$pqs$-systems. Then $\tau_1$ and $\tau_2$ are unitarily similar if
and only if their transfer functions coincide in some neighborhood
of zero.
\end{proposition}

%\section{The class ${\bf S}^{qs} $ and its properties}

\subsection{Transfer functions of $pqs$-systems}

Let $\tau=\{T;\sN,\sN,\sH\}$  be a $pqs$-system of the form
\eqref{abcd} and assume that $T$ is represented in the form
\eqref{MCE}. Then the transfer function $\Theta(\lambda)$ of $\tau$
has the form
\begin{equation}
\label{trfpqs} \Theta(\lambda)= K\Phi_A(\lambda)K^*+D_{K^*}XD_{K^*},
\quad \lambda \in \mathbb{D},
\end{equation}
where $\Phi_A(\lambda)$ is the characteristic function of the
selfadjoint contraction $A$; see \eqref{CHARFUNC1}. The function
$\Phi_A(\lambda)$ is holomorphic on $\dT\setminus\{-1,1\}$ and, in
fact, it belongs to Herglotz-Nevanlinna class on
$\Ext\{(-\infty,-1]\cup[1,\infty)\}$. Furthermore, $\Phi_A(\lambda)$
has nontangential strong limit values $\Phi_A(\pm 1)=\pm I_{\sD_A}$;
see e.g. \cite[Theorem~2.3]{AHS1}. Consequently, the limit value
$\Phi_A(\xi)$ is unitary for every $\xi\in \dT$ (see
\eqref{chardef}, \eqref{chardef*}), in particular, $\Phi_A(\lambda)$
is bi-inner. It follows from \eqref{trfpqs} that $\Theta(\lambda)$,
initially defined on $\mathbb{D}$, admits a holomorphic continuation
onto $\Ext\{(-\infty,-1]\cup[1,\infty)\}$. Furthermore,
$\Theta(\lambda)$ has nontangential strong limit values $\Theta (\pm
1)$ at $\pm 1$ which are given by
\begin{equation}\label{vaasa}
 \Theta(1)=KK^*+D_{K^*}XD_{K^*},\quad
 \Theta(-1)=-KK^*+D_{K^*}XD_{K^*}.
\end{equation}

Define the function $W(\lambda)$ by
\begin{equation}\label{W}
 W(\lambda)=\Theta(\lambda)-\Theta(0), \quad \lambda \in
 \Ext\{(-\infty,-1]\cup[1,\infty)\}.
\end{equation}
Since
\begin{equation}\label{W0}
 \Theta(0)=-KAK^*+D_{K^*}XD_{K^*},
\end{equation} it follows that
\begin{equation}\label{WW}
 W(\lambda)=\lambda K\left(I-\lambda A\right)^{-1}D^2_AK^*, \quad
 \lambda \in \Ext\{(-\infty,-1]\cup[1,\infty)\}.
\end{equation}
Hence $W^*(\ovl{\lambda})=W(\lambda)$ and
\[
\frac{W(\lambda)-W^*(\xi)}{\lambda-\ovl\xi}=\left\{
\begin{array}{l}
 KD_A(I-\lambda A)^{-1}(I-\ovl\xi A)^{-1}D_AK^*,\; \ovl\xi\ne \lambda,\\
 KD_A(I-\lambda A)^{-2}D_AK^*,\;\ovl\xi=\lambda.
\end{array}\right.
\]
Therefore $W(\lambda)$ is an operator-valued Herglotz-Nevanlinna
function with a holomorphic continuation onto
$\Ext\left\{(-\infty,-1]\cup [1,\infty)\right\}$. From \eqref{vaasa}
one sees that the strong limit values $W(\pm 1)$ exist and that they
are given by
\[
 W(1)=K(I+A)K^*, \quad W(-1)=-K(I-A)K^*.
\]
Hence,
\[
\frac{W(1)+W(-1)}{2}=KAK^*, \quad
I-\frac{W(1)-W(-1)}{2}=I-KK^*=D^2_{K^*}\ge 0.
\]
Since $X$ in \eqref{W0} is a contraction in $\sD_{K^*}$, these
identities show that
\begin{equation}\label{BW0}
\Theta(0)\in\bB\left(-\frac{W(1)+W(-1)}{2},\; I-\frac{W(1)-W(-1)}{2}
\right).
\end{equation}
Here $\bB(S,R)=\{\,S+R^{1/2}XR^{1/2}\in \bL(\sN):\, X \text{ a
contraction in }\, \bL(\cran R)\,\}$ stands for the \textit{operator
ball} with center $S\in \bL(\sN)$ and left and right radii $R\ge 0$.

\begin{definition}\label{DEFQS}
Let $\sN$ be a Hilbert space. The \textit{class} ${\bf S}^{qs}(\sN)$
consists of all $\bL(\sN)$-valued functions $\Theta(\lambda)$,
defined on $\dD$, such that
\begin{enumerate}%\def\labelenumi{\rm (\roman{enumi})}
\item[(S1)]
$W(\lambda)=\Theta(\lambda)-\Theta(0)$ is a Herglotz-Nevanlinna
function with a holomorphic continuation onto the domain
$\Ext\left\{(-\infty,-1]\cup [1,\infty)\right\}$;
\item[(S2)]
the strong limit values $W(\pm 1)$ exist and $W(1)-W(-1)\le 2I$;
\item[(S3)]
$\Theta(0)$ belongs to the operator ball in \eqref{BW0}.
\end{enumerate}
\end{definition}

The following proposition is now clear.

\begin{proposition}\label{recsys0}
Let $\tau=\{T;\sN,\sN,\sH\}$ be a $pqs$-system. Then its transfer
function $\Theta(\lambda)$ belongs to ${\bf S}^{qs}(\sN)$.
\end{proposition}

\section{The class ${\bf S}^{qs}$ and its realization via passive systems}

\subsection{The realization of the class ${\bf S}^{qs}$}

The next theorem is a converse to Proposition \ref{recsys0}. In its
proof a minimal $pqs$-system is constructed explicitly via an
operator representation of the Herglotz-Nevanlinna function
$W(\lambda)=\Theta(\lambda)-\Theta(0)$.

\begin{theorem}
\label{recsys} Let $\sN$ be a Hilbert space and let
$\Theta(\lambda)\in{\bf S}^{qs}(\sN)$. Then  $\Theta(\lambda)$  is
the transfer function of a minimal $pqs$-system
$\tau=\{T;\sN,\sN,\sH\}$.
\end{theorem}

\begin{proof}
Assume that $\Theta(\lambda)\in{\bf S}^{qs}(\sN)$. By the condition
(S1) the function
\[
\wt W(z):=-W(1/z),\quad z\in \Ext[-1,1],
\]
is a Herglotz-Nevanlinna function of the class
$\mathbf{N}_\sN[-1,1]$  with $\wt W(\infty)=0$, see \cite{AHS1}. It
follows from the condition (S2) that the strong limit values $\wt
W(\pm 1)$ exist. Then according to \cite[Theorem~2.3]{AHS1} there
exist a Hilbert space $\wt\sH$, a selfadjoint contraction $\wt A$ in
$\wt\sH$, and an operator $\wt G\in\bL(\sN,\sD_{\wt A})$, such that
\[
\wt W(z)=\wt G^*(\wt A- zI)^{-1}(I-\wt A^2)\wt G,
\]
see \cite{AHS1}. It follows that
\[
W(-1)=-\wt W(-1)=-\wt G^*(I-\wt A)\wt G, \quad W(1)=-\wt W(1)=\wt
G^*(I+\wt A)\wt G.
\]
Consequently,
\[
 \frac{W(1)+W(-1)}{2}=\wt G^*\wt A\wt G,
\quad
 I-\frac{W(1)-W(-1)}{2}=I-\wt G^*\wt G.
\]
The condition $W(1)-W(-1)\le 2I$ implies that $\wt G$ is
contractive. The condition (S3) means that $\Theta(0)=-\wt G^*\wt
A\wt G+D_{\wt G}\wt XD_{\wt G}$ for some contraction $\wt X$ in the
Hilbert space $\sD_{\wt G}$. Define in the Hilbert space
$\wt\cH=\sN\oplus \wt\sH$ the operator $\wt T$ by
\[
 \wt T=\begin{pmatrix} -\wt G^*\wt A\wt G+ D_{\wt G}\wt XD_{\wt G}&\wt G^*D_{\wt A}
     \cr D_{\wt
A}\wt G& \wt A
       \end{pmatrix}.
\]
Then $\wt T$ is a $qsc$-operator, $\ran(\wt T-\wt T^*)\subset \sN$,
and the operator $\wt T$ defines a $pqs$-system $\wt\tau =\{\wt
T;\sN,\sN,\wt\sH \}$; cf. Theorem \ref{MCE}.  The corresponding
transfer function is given by
\[
\Theta_{\wt \tau}(\lambda) =\wt G^*\left(-\wt A+\lambda
\left(I-\lambda \wt A\right)^{-1}
 D^2_{\wt A}\right)\wt G+D_{\wt G}\wt XD_{\wt G},
\quad  \lambda \in \dD.
\]
Therefore, $\Theta_{\wt
\tau}(\lambda)=\Theta(0)+W(\lambda)=\Theta(\lambda)$, $\lambda \in
\dD$. This means that the function $\Theta(\lambda)$ can be realized
as the transfer function of the $pqs$-system $\wt\tau$.  Finally,
replacing $\wt\tau$ by the system $\wt\tau^s$,  cf.
Proposition~\ref{simple1}, one obtains a minimal $pqs$-system.  The
corresponding transfer function still coincides with the function
$\Theta(\lambda)$.
\end{proof}

Observe that Theorem \ref{recsys} implies that the class ${\bf
S}^{qs}(\sN)$ is a subclass of the Schur class ${\bf S}(\sN)$.
Furthermore, the proof shows that a function $\Theta(\lambda)$ from
the class ${\bf S}^{qs}(\sN)$ admits the integral representation
\[
\Theta(\lambda)=\Theta(0)+\lambda\,\int_{-1}^{1}\frac{1-t^2}{1-t\lambda}\,d\Sigma(t),
\]
where $\Sigma(t)$ is a non-decreasing $\bL(\sN)$-valued function
with bounded variation, $\Sigma(-1)=0$, $\Sigma(1)\le I_\sN,$ and
\[
\left|\left(\left(\Theta(0)+\int_{-1}^{1}t\,d\Sigma(t)\right)f,g\right)\right|^2\le
\left(\left(I-\Sigma(1)\right)f,f\right)\,\left(\left(I-\Sigma(1)\right)g,g\right),
\quad f,g\in\sN.
\]

\begin{corollary}\label{IN}
Let $\sN$ be a Hilbert space and let $\Theta(\lambda) \in{\bf
S}^{qs}(\sN)$. If $\varphi_\Theta(\lambda)=0$
($\psi_\Theta(\lambda)=0$) then $\Theta(\lambda)$ is inner
(co-inner).
\end{corollary}

\begin{proof} By Theorem \ref{recsys} there exists a minimal
$pqs$-system $\tau=\{T,\sN,\sN,\sH\}$ with transfer function
$\Theta(\lambda)$. Now the statement follows from
Proposition~\ref{NorIN}.
\end{proof}

\begin{theorem}\label{INNER}
Let $\sN$ be a Hilbert space and let $\Theta(\lambda) \in{\bf
S}^{qs}(\sN)$. Then:
\begin{enumerate}\def\labelenumi{\rm (\roman{enumi})}
\item
if $\Theta(\lambda)$ is inner then% if and only if
\begin{equation}
\label{INNER1}
\begin{split}
&\left(\frac{\Theta(1)-\Theta(-1)}{2}\right)^2
=\frac{\Theta(1)-\Theta(-1)}{2},\\
&(\Theta(1)+\Theta(-1))^*(\Theta(1)+\Theta(-1))
=4I_\sN-2\left(\Theta(1)-\Theta(-1)\right);
\end{split}
\end{equation}
\item
if $\Theta(\lambda)$ is co-inner then% if and only if
\begin{equation}
\label{COINNER1}
\begin{split}
&\left(\frac{\Theta(1)-\Theta(-1)}{2}\right)^2
=\frac{\Theta(1)-\Theta(-1)}{2},\\
&(\Theta(1)+\Theta(-1))(\Theta(1)+\Theta(-1))^*
=4I_\sN-2\left(\Theta(1)-\Theta(-1)\right);
\end{split}
\end{equation}
\item
if \eqref{INNER1} (\eqref{COINNER1}) holds and $\Theta(\xi)$ is
isometric (co-isometric)  for some $\xi \in\dT$, $\xi\ne \pm 1$,
then $\Theta(\lambda)$ is inner (co-inner).
\end{enumerate}
\end{theorem}

\begin{proof}
Since $\Theta(\lambda) \in {\bf S}^{qs}(\sN)$, it is the transfer
function of a minimal $pqs$-system $\tau=\{T,\sN,\sN,\sH\}$. The
operator $T$, being quasi-selfadjoint, has the form \eqref{MCE} and
$\Theta(\lambda)$ is given by \eqref{trfpqs} with a holomorphic
continuation into the domain $\Ext\{(-\infty,-1]\cup [1,\infty)\}$.
Since $\Phi_A(\xi)$ is unitary for every $\xi\in\dT$, it follows
from \eqref{deffunc1} and \eqref{deffunc2} in Theorem~\ref{Schur}
and the definitions \eqref{fi} and \eqref{psi} that for all $h \in
\sN$ and $\xi\in\dT$
\begin{equation}
\label{deffunc3}
\begin{split}
&\left\|D_{\Theta (\xi)}h\right\|^2= ||D_{X}D_{K^*}h||^2
+\left\|\left(D_{K}\Phi_{A}(\xi)K^*-K^*XD_{K^*}
\right)h\right\|^2,\;  \\
&\left\|D_{\Theta^* }(\xi)h\right\|^2= ||D_{X^*}D_{K^*}h||^2
+\left\|\left(D_{K}\Phi_{A}(\ovl\xi)K^*-K^*X^*D_{K^*}
\right)h\right\|^2.
\end{split}
\end{equation}

(i) Suppose that $\Theta(\lambda)$ is inner. Then \eqref{deffunc3}
shows that
\[
\left\{
\begin{array}{l}
 D_XD_{K^*}=0,\\
 D_K\Phi_A(\xi)K^*=K^*XD_{K^*},\quad \xi\in\dT.
\end{array}
\right.
\]
The last equality yields that $D_{K}\Phi_A(\lambda)K^*=K^*XD_{K^*}$
for all $\lambda\in\dD$. Since
\[
\Phi_A(\lambda)=-A+\sum\limits_{n=0}^\infty \lambda^{n+1}
A^{n}D^2_A,
\]
it follows that $D_{K}D_AA^nD_AK^*=0$, $n \in \dN_0$. The minimality
of $\tau$ implies that $D_Kf=0$ for all $f\in\sD_A$; see
\eqref{neww}. Hence, $K$ is isometric and $D_{K^*}$ is an orthogonal
projector in the subspace $\sN$. Due to the identity $D_XD_{K^*}=0$,
$X\in\bL(\sD_{K^*})$ is isometric (here possibly $\sD_{K^*}=\{0\}$).
Now from the equalities in \eqref{vaasa} one obtains
\begin{equation}\label{UKSI}
 KK^*=\frac{\Theta(1)-\Theta(-1)}{2}
\end{equation}
and
\begin{equation}\label{KAKSI}
 XD_{K^*}=D_{K^*}XD_{K^*}=\frac{\Theta(1)+\Theta(-1)}{2}.
\end{equation}
%since $D_{K^*}$ is an orthogonal projector and $X\in\bL(\sD_{K^*})$.
These identities %\eqref{UKSI}, \eqref{KAKSI}, and
together with the equality $X^*X=I_{\sD_{K^*}}$ lead to
\eqref{INNER1}.

(ii) The proof is similar to that of (i).

(iii) Assume that \eqref{INNER1} holds and that $\Theta(\xi)$ is
isometric for some $\xi\in\dT$, $\xi\ne \pm 1$. Due to
\eqref{INNER1} $D_{K^*}$ is an orthogonal projector in $\sN$ and $X$
is isometric in $\sD_{K^*}$. Hence, $K^*\in\bL(\sN,\sH)$ is a
partial isometry and, moreover, $D_KK^*=K^*D_{K^*}=0$ and
$K^*XD_{K^*}=0$. Since $\Theta(\xi)$ is isometric, \eqref{deffunc3}
gives $D_K\Phi_A(\xi)K^*=0$. Furthermore, since
\[
 -A+\xi(I_\sH-A^2)(I_\sH-\xi A)^{-1}
 =\ovl\xi+(\xi-\ovl\xi)(I_\sH-\xi A)^{-1},
\]
the equality $D_K\Phi_A(\xi)K^*=0$ with $\xi\neq\ovl\xi$ implies
that $D_K(I_\sH-\xi A)^{-1}K^*=0$, i.e.
\[
 (I_\sH-\xi A)^{-1}\ran K^*\subset \ran K^*.
\]
Thus $A(\ran K^*)\subset \ran K^*$ (since $A=A^*$) and hence
$D_A(\ran K^*)\subset\ran K^*$, so that
\begin{equation}\label{innerpqs}
\cspan\{A^nD_{A}K^*\sN :\, n \in \dN_0 \}\subset \ran K^*.
\end{equation}
Consequently $\ran K^*=\sH$, i.e., $K\in\bL(\sH,\sN)$ is isometric.
Therefore, $D_K\Phi_A(\zeta)K^*=0$ for all $\zeta\in\dT$ and
$\Theta(\lambda)$ is inner in view of \eqref{deffunc3}. Similarly,
if \eqref{COINNER1} holds and $\Theta(\xi)$ is co-isometric for some
$\xi \in\dT$, $\xi\ne \pm 1$, then $\Theta(\lambda)$ is co-inner.
\end{proof}

\begin{theorem}\label{CONSREAL}
Let the Hilbert space $\sN$ be finite-dimensional and let
$\Theta(\lambda)$ be a nonconstant inner function from  ${\bf
S}^{qs}(\sN)$. Then $\Theta(\lambda)$ is rational and
\[
\Theta(\lambda)=\diag \left( \frac{\lambda-a_1}{1-\lambda\,a_1},
\frac{\lambda-a_2}{1-\lambda\,a_2}, \ldots,
\frac{\lambda-a_m}{1-\lambda\,a_m}, X \right)
\]
relative to some orthonormal basis in $\sN$. Here the not
necessarily distinct numbers $a_1,a_2,\ldots,a_m$ belong to
$(-1,1)$, and $X$ is a constant unitary matrix.
\end{theorem}

\begin{proof}
Since $\Theta(\lambda)\in{\bf S}^{qs}(\sN)$, it is the transfer
function of a minimal $pqs$-system $\tau=\{T,\sN,\sN,\sH\}$. As $T$
is quasi-selfadjoint, it has the form \eqref{MCE}, and
$\Theta(\lambda)$ is given by \eqref{trfpqs}. Since $\sN$ is
finite-dimensional and $\Theta(\lambda)$ is inner, $\Theta(\lambda)$
is automatically bi-inner. Then in \eqref{einz} one has
$\f_\Theta(\lambda)=0$ and $\psi_\Theta(\lambda)=0$. Thus by
\cite[Theorem~1.1]{AHS2} $\tau$ is conservative (in fact, this
conclusion can be derived also from the proof of Theorem~\ref{INNER}
above by applying Corollary~\ref{Corvaasa}). As $\Theta(\lambda)$ is
nonconstant, $A$ is non-isometric, $K$ is isometric, and $X$
appearing in \eqref{MCE} is unitary in $\sD_{K^*}$. Since $\tau$ is
minimal, Proposition \ref{simple1} shows that ${\sD}_A=\sH=\sH^s$.
Hence $K\in\bL(\sD_A,\sN)$ isometric implies that $\dim \sH\le
\dim\sN<\infty$, so that $\Theta(\lambda)$ is rational.  Let
$\sN_0=\ran K$, $\sN_1=\sN\ominus\sN_0=\ker K^*.$ Then
$K^*\sN_0=\sH$ and
\[
\Theta(\lambda)\uphar\sN_1=X:\sN_1\to\sN_1,\;
\Theta(\lambda)\uphar\sN_0=K\Phi_A(\lambda)K^*\uphar\sN_0:\sN_0\to\sN_0.
\]
Suppose that $\dim\sH=m$ and  that $a_1,\ldots, a_m$ are the
eigenvalues of $A$. Choose an orthonormal basis in $\sH$ consisting
of eigenvectors of $A$. Then $K$ maps this basis onto some
orthonormal basis in $\sN_0$. With respect to this basis the matrix
$\Theta(\lambda)\uphar\sN_0$ is diagonal with entries
$\theta_{kk}(\lambda)$, $k=1,\ldots, m$, and since $\Phi_A(\lambda)$
has the form \eqref{CHARFUNC1}, it follows that
\[
\Theta_{kk}(\lambda)=-a_k+\frac{(1-a^2_k)\lambda}{1-\lambda\,
a_k}=\frac{\lambda-a_k}{1-\lambda\,a_k}.
\]
This completes the proof.
\end{proof}

\subsection{Bi-inner dilations of  functions from the class ${\bf S}^{qs}$}

The function $\Theta(\lambda) \in {\bf S}(\sM,\sN)$  is said to have
an \textit{inner dilation} if there exists a function
$\Theta_r(\lambda)$  such that
\[
{\bf
\Theta}(\lambda)=\begin{pmatrix}\Theta(\lambda)\cr\Theta_r(\lambda)
\end{pmatrix}\in {\bf S}(\sM,\sN\oplus\sL)
\]
is inner. The function $\Theta(\lambda) \in {\bf S}(\sM,\sN)$  is
said to have a \textit{co-inner dilation} if there exists a function
$\Theta_l(\lambda)$ such that
\[
{\bf
\Theta}(\lambda)=\begin{pmatrix}\Theta(\lambda)&\Theta_l(\lambda)
\end{pmatrix}\in {\bf S}(\sM\oplus\sK,\sN)
\]
is co-inner. The function $\Theta(\lambda) \in {\bf S}(\sM,\sN)$  is
said to have a \textit{bi-inner dilation} if there exist functions
$\Theta_{11}(\lambda)$, $\Theta_{22}(\lambda)$,  and
$\Theta_{21}(\lambda)$ such that
\begin{equation}\label{bi}
{\bf \Theta}(\lambda)
=\begin{pmatrix}\Theta(\lambda)&\Theta_{12}(\lambda)&\cr\Theta_{21}(\lambda)
&\Theta_{22}(\lambda)\end{pmatrix}\in {\bf
S}(\sM\oplus\sK,\sN\oplus\sL)
\end{equation}
is bi-inner. Recall the following result due to Arov \cite{Arov};
cf. \cite{ArSt}.

\begin{proposition}\label{biinner}
Let $\tau=\{T;\sM,\sN,\sH\}$ be a passive system with transfer
function $\Theta (\lambda)$. Then:
\begin{enumerate}\def\labelenumi{\rm (\roman{enumi})}
\item
if $\tau$ is strongly stable, then $\Theta(\lambda)$ has an inner
dilation;
\item
if $\tau$ is strongly co-stable, then $\Theta(\lambda)$ has a
co-inner dilation;
\item
if $\tau$ is strongly stable and strongly co-stable, then
$\Theta(\lambda)$ has a bi-inner dilation.
\end{enumerate}
\end{proposition}

In \cite{AHS2} this result was proved using the parametrization in
Theorem \ref{ParContr} and   \eqref{exptran}.  Since a function
$\Theta(\lambda) \in {\bf S}^{qs}(\sN)$ can be realized as the
transfer function of a minimal $pqs$-system, it is strongly stable
and strongly co-stable by Corollary \ref{Norcoro}. Hence it admits a
bi-inner dilation by Proposition \ref{biinner}.

\begin{proposition}\label{BIIN}
Among the bi-inner dilations of $\Theta(\lambda) \in{\bf
S}^{qs}(\sN)$ there exists a bi-inner dilation from the class ${\bf
S}^{qs}(\sN)$.
\end{proposition}

\begin{proof}
Since $\Theta(\lambda) \in {\bf S}^{qs}(\sN)$, it is the transfer
function of a minimal $pqs$-system $\tau=\{T,\sN,\sN,\sH\}$ by
Theorem~\ref{recsys}. The operator $T$, being quasi-selfadjoint, has
the form \eqref{MCE} and therefore $\Theta(\lambda)$ is given by
\eqref{trfpqs}. Define the following functions
\[
\Theta_{21}(\lambda):=\begin{pmatrix}
D_{K}\Phi_{A}(\lambda)K^*-K^*XD_{K^*} \cr
-D_XD_{K^*}\end{pmatrix}:\sN\to\begin{pmatrix}\sD_{K}\\
\sD_{K^*}\end{pmatrix},\quad \lambda\in\dD,
\]
\[
 \Theta_{12}(\lambda):=\begin{pmatrix}  K\Phi_{A}(\lambda)
D_{K}-D_{K^*}XK &D_{K^*}D_{X^*}\end{pmatrix}:
\begin{pmatrix}\sD_{K}\\ \sD_{K^*}\end{pmatrix}\to\sN,
\quad \lambda\in\dD,
\]
\[
\Theta_{22}(\lambda)
 =\begin{pmatrix} K^*XK+D_{K}\Phi_{A}(\lambda)D_{K}& -K^*D_{X^*}
  \cr D_XK&X^*\end{pmatrix}:\begin{pmatrix}\sD_{K}\\ \sD_{K^*}\end{pmatrix}\to
  \begin{pmatrix}\sD_{K}\\ \sD_{K^*}\end{pmatrix},\quad \lambda \in\dD.
  \]
Let $\sK=\sL=\sD_{K}\oplus\sD_{K^*}$ and $\sV=\sN\oplus\sL$, and let
${\bf \Theta}(\lambda)$ be defined by \eqref{bi}. Furthermore,
define the operator $\bf T$ by
\[
{\bf T}=\begin{pmatrix}{\bf D}& {\bf C}\cr{\bf B}&{\bf A}
\end{pmatrix}:\begin{pmatrix}{\sV}\\ \sH\end{pmatrix}\to
 \begin{pmatrix}{\sV}\\ \sH\end{pmatrix},
\]
where ${\bf A}=A$, ${\bf D}={\bf \Theta}(0)$, and
\[
{\bf B}=D_A\begin{pmatrix}K^*&D_K&0\end{pmatrix}:\begin{pmatrix}\sN\\
\sD_K\\ \sD_{K^*}\end{pmatrix}\to\sH, \,\, {\bf C}={\bf
B}^*=\begin{pmatrix}K\cr D_K\cr 0\end{pmatrix}D_A:\sH\to
\begin{pmatrix}\sN\\ \sD_K\\ \sD_{K^*}\end{pmatrix}.
\]
Simple calculations show that the operator ${\bf T}$ is unitary and
quasi-selfadjoint. Hence, the system $\eta=\left\{{\bf
T};{\sV},{\sV},\sH\right\}$ is conservative. Since $A\in C_{00}$,
the system $\eta$ is minimal; see Corollary~\ref{stablecor} and e.g.
\cite[Proposition~5.2]{AHS2}. In addition it is easy to see that the
transfer function of $\eta$ coincides with ${\bf\Theta}(\lambda)$
and therefore ${\bf\Theta}(\lambda) \in {\bf S}^{qs}({\sV})$ by
Proposition~\ref{recsys0}. Since $\eta$ is minimal and conservative,
${\bf \Theta}(\lambda) \in {\bf S}({\sV})$ is a bi-inner dilation of
${\Theta}(\lambda) \in {\bf S}({\sN})$ in view of
\cite[Corollary~5.3]{AHS2}.
\end{proof}

\begin{remark}
With straightforward calculations it is easy to see directly that
${\bf \Theta}(\lambda) \in {\bf S}({\sV})$ as defined explicitly in
Proposition~\ref{BIIN} is bi-inner; cf.
\cite[Proposition~7.1]{AHS2}. Furthermore, using the explicit
formula of ${\bf \Theta}(\lambda)$ one can also calculate the
function ${\bf W}(\lambda):={\bf \Theta}(\lambda)-{\bf \Theta}(0)$
as introduced in Definition~\ref{DEFQS}. Observe, that
\[
{\bf W}(\lambda): \begin{pmatrix}\sN\\ \sD_K\\
\sD_{K^*}\end{pmatrix}\to
\begin{pmatrix}\sN\\ \sD_K\\ \sD_{K^*}\end{pmatrix}
\]
and that
\[
{\bf\Theta}(0)=
\begin{pmatrix}-KAK^*+D_{K^*}XD_{K^*}&-KAD_K-D_{K^*}XK&D_{K^*}D_{X^*}\cr
-D_KAK^*-K^*XD_{K^*}&K^*XK-D_KAD_K&-K^*D_{X^*}\cr
-D_XD_{K^*}&D_XK&X^*.
\end{pmatrix}.
%:\begin{pmatrix}\sN\\ \sD_K\\ \sD_{K^*}\end{pmatrix}\to
%\begin{pmatrix}\sN\\ \sD_K\\ \sD_{K^*}\end{pmatrix},
\]
Thus,
\[
{\bf W}(\lambda) =\begin{pmatrix} \lambda KD^2_{A}(I_\sH-\lambda
A)^{-1}K^*& \lambda KD^2_{A}(I_\sH-\lambda A)^{-1}D_K&0\cr \lambda
D_KD^2_{A}(I_\sH-\lambda A)^{-1}K^*&\lambda D_KD^2_{A}(I_\sH-\lambda
A)^{-1}D_K&0\cr 0&0& 0
\end{pmatrix}
\]
which implies that the Nevanlinna kernel
\[
 \cfrac{{\bf W}(\lambda)-{\bf W}^*(\lambda)}{\lambda-\bar \lambda}
\]
is given by
\[
\begin{small}
\begin{pmatrix}  KD^2_{A}(I_\sH-\bar\lambda A)^{-1}(I_\sH-\lambda A)^{-1}K^*&
 KD^2_{A}(I_\sH-\bar\lambda A)^{-1}(I_\sH-\lambda A)^{-1}D_K&0\cr
 D_KD^2_{A}(I_\sH-\bar\lambda A)^{-1}(I_\sH-\lambda A)^{-1}K^*&
 D_KD^2_{A}(I_\sH-\bar\lambda A)^{-1}(I_\sH-\lambda A)^{-1}D_K&0\cr
0&0&0\end{pmatrix}.
\end{small}
\]
Therefore, ${\bf W}(\lambda)$ is a Herglotz-Nevanlinna function
defined on $\Ext\{(-\infty,-1]\cup[1,\infty)\}$, which reflects the
defining properties of the class ${\bf S}^{qs}(\sN)$ in
Definition~\ref{DEFQS}.
\end{remark}

\section{Minimal systems with transfer functions of the class ${\bf S}^{qs}$}

In this section the class of minimal systems
$\dot\tau=\{T,\sN,\sN,\sH\}$ of the form \eqref{passive} are
considered. The system  $\dot\tau$ is not assumed to be passive, so
that the bounded operator $T$ in \eqref{abcd0} need not be
contractive. The next result can be seen as an extension of the
unitary similarity result for $pqs$-systems in Proposition~\ref{US}.

\begin{theorem}\label{opt}
Let $\Theta(\lambda)\in {\bf S}^{qs}(\sN)$ and let
\begin{equation}\label{dot}
\dot{\tau}
=\left\{\begin{pmatrix}\dot{D}&\dot{C}\cr\dot{B}&\dot{A}\end{pmatrix};\sN,\sN,\sH\right\}.
\end{equation}
be a minimal, not necessarily passive, system whose transfer
function coincides with $\Theta(\lambda)$ in some neighborhood of
zero. Then there exists a positive selfadjoint operator $S$ in $\sH$
such that
\[
 \dot{A}^*S=S\dot{A},\quad \dot{C}^*=S\dot{B},
\]
and such that
\[
\wh \tau_0=\left\{\,\begin{pmatrix}D&\dot{C}S^{-1/2}\cr
S^{1/2}\dot{B} &S^{1/2}\dot{A} S^{-1/2}
\end{pmatrix};\sN,\sN,\sH\,\right\}
\]
is a minimal $pqs$-system with transfer function $\Theta(\lambda)$.
Moreover,   $\dot{\tau}$ becomes a $pqs$-system with respect to the
inner product in $\dom S^{1/2}\subset\sH$ given by
\[
(\varphi,\psi)_{S^{1/2}}=(S^{1/2}\varphi,S^{1/2}\psi)_{\sH}, \quad
\varphi, \psi \in \dom S^{1/2}.
\]
Furthermore, if   $\Theta(\lambda)$ is inner (co-inner), then all
minimal passive realizations of  $\Theta(\lambda)$ are unitarily
similar isometric (co-isometric) $pqs$-systems.
\end{theorem}

\begin{proof}
Since $\Theta(\lambda)\in {\bf S}^{qs}(\sN)$, it is the transfer
function of a minimal $pqs$-system $\tau_0$ of the form
\begin{equation}\label{tau1}
\tau_0=\left\{\begin{pmatrix}D&C\cr B&A\end{pmatrix}; \sN,\sN,
{\stackrel{\rm o}\sH}\right\}.
\end{equation}
Here  $ A$ is a selfadjoint contraction in the Hilbert space
${\stackrel{\rm o} {\sH}}$ and $C=B^*$. Furthermore,
$\Theta(\lambda)$ has the form
\[
\Theta(\lambda)=W(\lambda)+D,
\]
where $W(\lambda)=\lambda B^*(I-\lambda A)^{-1}B$, $\lambda \in
\dD$, is a Herglotz-Nevanlinna function on the domain
$\Ext\left\{(-\infty,-1]\cup [1,\infty)\right\}$ and $W(0)=0$. Since
the transfer function of $\dot{\tau}$ in \eqref{dot} is of the form
\[
\Theta(\lambda)=D+\lambda\dot{C}(I-\lambda \dot{A})^{-1}\dot B,
\]
one concludes that $\dot{D}=D$ and in some neighborhood of zero
\begin{equation}\label{B}
 \dot{C}(I-\lambda \dot{A})^{-1}\dot{B} = B^*(I-\lambda A)^{-1}B.
\end{equation}
Consequently,
\begin{equation}
\label{nevtr} \dot{C}\dot{A}^k\dot{B}
 = B^*A^k B,\quad k \in \dN_0.
\end{equation}
%\marginpar{I stopped here with section 5}STILL TO USE The system
%$\dot{\tau}$ is controllable and observable, i.e.
%\[
%\sH=\cspan\{\ran \dot{A}^n\dot{B}:\,n\in \dN_0  \}= \cspan\{\ran
%\dot{A}^{*n}\dot{C}^*:\, n \in \dN_0 \}.
%\]
Define the linear operator $Y$ by
\begin{equation}
\label{Y}
 Y=\left\{\, \left\{\sum\limits_{k=0}^n \dot{A}^k\dot{B}u_k,
 \sum\limits_{k=0}^n A^k Bu_k\right\}:\,u_0, u_1,..., u_n\subset \sN,
 \quad n \in \dN_0 \,\right\}.
\end{equation}
%\begin{equation}
%\label{Y} Y\left(\sum\limits_{k=0}^n \dot{A}^k\dot{B}u_k\right)=
%\sum\limits_{k=0}^n A^k Bu_k,
%\end{equation}
It follows from \eqref{degn0} that
\begin{equation}\label{Y*}
 \left\{\,\sum\limits_{k=0}^m A^k Bv_k,
\sum\limits_{k=0}^m \dot{A}^{*k}\dot{C}^*v_k\,\right\} \in Y^*,
 \quad v_1,\ldots,v_m \in \sN, \quad m \in \dN_0.
\end{equation}
By minimality of $\dot{\tau}$ and $\tau_0$, $Y$ and its adjoint
$Y^*$ have a dense domain and dense range. In particular, the
operator $Y$ is closable and the closure $\ovl Y=Y^{**}$ is a
densely defined operator with dense range. Definition~\ref{Y} and
the relations \eqref{nevtr}, \eqref{Y*} yield the equalities
%\begin{equation}
%\label{simx} \left\{
%\begin{split}
%&{\dot{A}^*}\ovl Xv=\ovl X Av,\; v\in\dom\ovl X,\\
%&{\dot{B}}^*\ovl X v=B^* v,\; v\in\dom\ovl X,\\
%&XB={\dot{C}}^*,
%\end{split}
%\right.
%\end{equation}
\begin{equation}
\label{simy} \left\{
\begin{array}{l}
 A\ovl Yu=\ovl Y{\dot{A}}u,\; u\in\dom\ovl Y,\\
 \dot{C} u =B^*\ovl Yu,\; u\in\dom\ovl Y,\\
 Y\dot{B}=B,
\end{array}
\right.
\end{equation}
which means that the systems $\dot{\tau}$ and $\tau_0$ are weakly
similar; see \eqref{pseudosim}, cf. also \cite{A}. Now define a
positive selfadjoint operator $S$ by $S=Y^*\ovl Y$. Since
$Y^*B={\dot{C}}^*$ and $Y\dot{B}=B$, one obtains
$S\dot{B}={\dot{C}}^*$. Represent $\ovl Y$ in the form $\ovl
Y=US^{1/2}$, where $U:\sH\to \stackrel{0}{\sH} $ is an isometry and
$\dom S^{1/2}=\dom\ovl Y$. Because $\cran Y=\stackrel{0}{\sH}$, the
operator $U$ is a unitary from $\sH$ onto $\stackrel{0}{\sH}$.
Define a selfadjoint contraction $\wh A$ in $\sH$ by the equality:
\[
\wh A=U^{-1}A U.
\]
Then the equality $\ovl Y{\dot{A}}=A\ovl Y$ gives
\[
 S^{1/2 }{\dot {A}}v =\wh AS^{1/2}v,\; v\in\dom S^{1/2},
 \]
or, equivalently,
\[
\wh Av=S^{1/2}\dot{A}S^{-1/2}v,\; v\in\dom S^{-1/2}.
\]
Let $\wh B=U^{-1}B.$ From $Y{\dot{B}}=B$ one obtains $S^{1/2}{\dot
{B}}=\wh B$ and the relation $S\dot{B}={\dot{C}}^*$ yields
$S^{-1/2}{\dot{C}}^*=\wh B.$ Therefore $\wh
B^*v={\dot{C}}S^{-1/2}v,$ $v\in\dom S^{-1/2}.$ The system
\[
\wh\tau_0=\left\{\begin{pmatrix}D&\wh B^*\cr \wh B& \wh
A\end{pmatrix};\sN,\sN,\sH\right\}
\]
is unitarily equivalent to $\tau_0$. Hence $\wh\tau_0$ is a minimal
$pqs$-system with the transfer function $\Theta(\lambda).$ Let
$v\in\dom S,$ $u\in\dom S^{1/2}$. Then
\begin{equation}\label{SAeq}
\begin{split}
 \left(S^{1/2}{\dot{A}}v, S^{1/2}u\right)
 &=\left(\wh AS^{1/2}v,S^{1/2}u\right)=\left(S^{1/2}v,\wh A S^{1/2}u\right)\\
 &=\left(S^{1/2}v,S^{1/2}{\dot{A}}u\right)=\left({\dot{A}}^*Sv,u\right).
\end{split}
\end{equation}
It follows that the vector $S^{1/2}{\dot{A}}v$ belongs to $\dom
S^{1/2}$ and
\[
S{\dot{A}}v={\dot{A}}^*Sv,\; v\in\dom S.
\]
%Suppose that the system $\dot\tau$ is optimal \cite{Arov},
%\cite{ArKaaP}, \cite{ArKaaP3} and minimal. It follows that the
%system
%$$\dot{\tau}^*=\left\{\begin{pmatrix}D^*&\dot{C}^*\cr \dot{B}^*&\dot{A}^*\end{pmatrix};\sN,\sN,\sH\right\}$$
%is the $(*)$ - optimal and minimal realization of the function
%$\Theta^*(\ovl\lambda)$ \cite{Arov}, \cite{ArKaaP}, \cite{ArKaaP3}.
%But the $pqs$-system
%\[
%\tau_0^*=\left\{\begin{pmatrix}D^*&B^*\cr
%B&A^*\end{pmatrix};\sN,\sN,\stackrel{0}{\sH}\right\}
%\]
%is minimal and has the transfer function $\Theta^*(\ovl\lambda)$. By
%definition of optimal and $(*)$ - optimal systems \cite{Arov},
%\cite{ArKaaP}, \cite{ArKaaP3} we obtain
%\[
%%\begin{equation}
%%\label{INEQ}
%\left\|\sum\limits_{k=0}^n \dot{A}^k\dot{B}u_k\right\|\le
%\left\|\sum\limits_{k=0}^n A^k Bu_k\right\|\le
%\left\|\sum\limits_{k=0}^n \dot{A}^{*k}\dot{C}^*u_k\right\|
%%\end{equation}
%\]
%for all $n=0,1,...$ and every choice of vectors $u_0,u_1,...,
%u_n\in\sN.$ It follows that the operators $ X^{-1}$, $Y^{-1}$, and
%$S^{-1}$ are contractions.
%\marginpar{left out too much?}
Consider $\dom S^{1/2}$ as a pre-Hilbert space equipped with the
inner product $(f,g)_{S^{1/2}}:=(S^{1/2}f,S^{1/2}g)$ and let
${\sH}_{S^{1/2}}$ be the completion of $\dom S^{1/2}$ with respect
to this inner product. Since $S^{1/2}\dot{A}v=\wh AS^{1/2}v$,
$v\in\dom S^{1/2}$, one obtains from \eqref{SAeq}
\[
(\dot{A}u,v)_{S^{1/2}}=(S^{1/2}\dot{A}v, S^{1/2}u)
 =(S^{1/2}v,S^{1/2}\dot{A}u)=(v,\dot{A}u)_{S^{1/2}},
\]
for all $v,u\in\dom S^{1/2}$. Thus, the operator $\dot{A}$ is
symmetric with respect to the inner product
$(\cdot,\cdot)_{S^{1/2}}$. Furthermore, the operator $\dot T$
defined via the block formula
\[ \dot{T} =\begin{pmatrix}D & \dot{C}\cr
\dot{B}& \dot{A}\end{pmatrix} :
\begin{pmatrix}\sN\\ {\sH}_{S^{1/2}}\end{pmatrix} \to
\begin{pmatrix}\sN\\{\sH}_{S^{1/2}}\end{pmatrix}
\]
is a $qsc$-operator. In fact, since $\dot{C}^*=S^{1/2}\wh B$ and
$S^{1/2}\dot B=\wh B$, the operator $\dot B:\sN\to\sH_{S^{1/2}}$ is
the adjoint of the operator $\dot C:\dom
S^{1/2}\subset\sH_{S^{1/2}}\to\sN$. Because the operator matrix
\[
\wh T =\begin{pmatrix}D&{\wh B}^*\cr \wh B & \wh A\end{pmatrix}:
\begin{pmatrix}\sN\\ {\sH}\end{pmatrix} \to
\begin{pmatrix}\sN\\ {\sH}\end{pmatrix}
\]
is a contraction, one obtains with $f\in\dom S^{1/2}$ and $u\in\sN$
\[
\begin{split}
&\left\|\dot{A}f+\dot B u\right\|^{2}_{S^{1/2}}+
\left\|\dot C f+Du\right\|^2 \\
&=\left\|S^{1/2}(\dot A f+\dot B u) \right\|^2+
\left\|\dot C f+Du\right\|^2 \\
&=\left\|\wh AS^{1/2}f+\wh Bu\right\|^2+
\left\|\wh B^*S^{1/2}f+Du\right\|^2 \\
&\le ||S^{1/2}f||^2+||u||^2=\left\|f\right\|^2_{S^{1/2}}+||u||^2.
\end{split}
\]
Therefore, after a renormalization of the state space by means of
the operator $S$ the system $\dot\tau$ becomes a $pqs$-system.

Finally, if the function $\Theta(\lambda) \in {\bf S}^{qs}(\sN)$ is
inner (co-inner) then $\varphi_\Theta(\xi)=0$ ($\psi_\Theta(\xi)=0$)
for almost all $\xi\in \dT$, see \eqref{einz}, and now it follows
from \cite[Theorem~1.1]{AHS2} that every two minimal passive
realizations of the inner (co-inner) function $\Theta(\lambda) \in
{\bf S}^{qs}(\sN)$ are unitarily similar isometric (co-isometric)
$pqs$-systems.
\end{proof}

\section{The class ${\bf S}^{qs}$ and $Q$-functions of
quasi-selfadjoint contractions} \label{sec6}

\subsection{Hermitian contractions}

Let $A_0$ be a Hermitian contraction in the Hilbert space $\cH$ with
$\dom A_0=\sH$ and let $\sN=\cH\ominus\sH$. Then $A:=P_\sH
A_0\uphar\sH$ is   selfadjoint  in the Hilbert space $\sH$. It
follows that there is a contraction $K\in\bL(\sD_{A},\sN)$ such that
\begin{equation}\label{anul}
A_0=A+KD_A.
\end{equation}
The operator  $A_0$ is said to be \textit{simple} if there is no
nonzero subspace in $\dom A_0$ which is invariant under $A_0$; cf.
\cite{KrO}. Since $A_0$ is Hermitian, simplicity of $A_0$ is
equivalent to $A_0$ being completely nonselfadjoint, i.e., $A_0$ has
no selfadjoint part.  Note that $A_0$ is simple if and only if the
subspace $\sH^s_0$ in \eqref{SP} coincides with $\sH$; cf.
\cite[Lemma~3.2]{AHS1}. For a Hermitian contraction $A_0=A+KD_{A}$
the restriction $A_0\uphar\sH^s_0$ is called the \textit{simple part
of $A_0$}.

A $qsc$-operator $T$ is said to be a \textit{quasi-selfadjoint
contractive extension} or \textit{$qsc$-extension} of a Hermitian
contraction $A_0$ if
\[
  A_0\subset T \mbox{ and } A_0\subset T^*,
\]
or, equivalently, if $\dom A_0\subset\ker(T-T^*)$, cf. \cite{ArTs1},
\cite{ArTs2}. A $qsc$-operator $T$ has always Hermitian restrictions
$A_0$ for which $T$ is a $qsc$-extension. Namely, with a subspace
$\sN\supset\ran(T-T^*)$ define
\[
\dom A_0=\sH:=\cH\ominus \sN,\quad A_0=T\uphar \dom A_0.
\]
Then $\dom A_0\subset\ker(T-T^*)$.

With $A_0$ and $K\in\bL(\sD_A,\sN)$ as in \eqref{anul}, the formula
\eqref{MCE} provides a one-to one correspondence between all
contractions $X\in\bL(D_{K^*})$ and all $qsc$-extensions of $A$. The
operator form of all $qsc$-extensions with their resolvents was
obtained in \cite{ArTs3}, \cite{ArTs2}. Clearly, the subspaces
$\sH'$ and $\sH''=\sH \ominus \sH'$, where
\begin{equation}
\label{HH}  \sH':=\cspan\left\{(T-zI)^{-1}\sN:\,|z|>1\right\}
      =\cspan\left\{T^n\sN :\,n \in \dN_0 \right\},
\end{equation}
are invariant with respect to $T$ and $T^*$, respectively. The
inclusion $\sN\subset \sH'$ implies that $\sH''\subset\sN^\perp=\dom
A\subset\ker(T-T^*)$. Therefore the restriction of $T^*$ to $\sH''$
is a selfadjoint operator in $\sH''$. The restriction $T\uphar
\sH'\,(=P_{\sH'}T\uphar\sH')$ is called the \textit{$\sN$-minimal
part of $T$}. Moreover, $T$ is said to be \textit{$\sN$-minimal} if
the equality $\sH=\sH'$ holds \cite{AHS1}. The subspaces $\sH'$ and
$\sH^s$ of $\cH=\sN\oplus\sH_0$ defined in \eqref{HH} and
\eqref{SP}, respectively, are connected by $ \sH'=\sN\oplus\sH^s$.
Every $qsc$-extension $T$ of the Hermitian contraction $A_0$ in
$\cH$ with $\dom A_0=\sH$ generates a $pqs$-system in the following
manner: let \eqref{abcd} be the block-operator representation of
$T$, then the system $\tau=\{T;\sN,\sN, \sH\}$ is a $pqs$-system.

\begin{proposition}\label{mpss}
Let $A_0$ be a Hermitian contraction in $\cH=\sN\oplus\sH$ with
$\sH=\dom A_0$, let $T$ be a $qsc$-extension of $A$ in $\cH$, and
let $\tau$ be the $pqs$-system generated by $T$ with the state space
$\sH$ and the input and the output space $\sN$. Then the following
statements are equivalent:
\begin{enumerate}\def\labelenumi{\rm (\roman{enumi})}
\item the Hermitian contraction $A_0$ is simple;
\item the $qsc$-extension $T$ of $A_0$ is $\sN$-minimal;
\item the $pqs$-system $\tau$ is minimal.
\end{enumerate}
\end{proposition}

\begin{proof} The equivalence (i) $\Leftrightarrow$ (ii) was proved in
\cite{AHS1}. For the proof of (i) $\Leftrightarrow$ (iii), observe
that the simplicity of $A_0$ is equivalent to $\sH_0^s=\sH$. Since
$\ran K^*\subset \sD_{A_0}$ and $\sD_{A_0}$ is invariant under
$A_0$, the inclusion $\sH_0^s\subset \sD_{A_0}$ holds. Hence,
$\sH_0^s=\sH$ implies in fact that $\ker D_{A_0}=\{0\}$, i.e.,
$\|Af\|<\|f\|$ holds for all $f\in\sH\setminus\{0\}$. Now the
assertion follows from Proposition~\ref{simple1}.
\end{proof}

\subsection{Transfer functions and $Q$-functions}

Let $T$ be a $qsc$-operator in a separable Hilbert space $\cH$ and
let $\sN$ be a subspace of $\sH$ such that $\sN\supset\ran(T-T^*)$.
The operator-valued function
\begin{equation}
\label{QF} Q_T(z)=P_\sN(T-zI)^{-1}\uphar \sN, \quad |z|>1,
\end{equation}
is said to be a \textit{$Q$-function} of $T$, cf. \cite{AHS1}.
Analytical properties of the $Q$-function of $qsc$-operators and its
applications to the parametrization of the resolvents of all
$qsc$-extenions of corresponding Hermitian contraction were
established in \cite{AHS1}.

Let $\sN$ be a Hilbert space. An operator-valued function $Q(z)$
with values in $\bL(\sN)$ and holomorphic outside $\dD$ is said to
belong to the class $\bQ(\sN)$ (see \cite[Section~6]{AHS1}) if
\begin{enumerate}
\item[(S1)]  $Q(z)$ has the asymptotic expansion
\begin{equation}
\label{EXP}
Q(z)=-\frac{1}{z}I+\frac{1}{z^2}F+o\left(\frac{1}{z^2}\right),\;
z\to \infty;
\end{equation}
\item[(S2)]  the  kernel
\[
%\fG(z,\xi)=
\frac{Q(z)-Q^*(\xi)-Q^*(\xi)(F-F^*)Q(z)}{z-\bar\xi}
%\quad |z|,|\xi|>1,
\]
is  nonnegative;
\item[(S3)] the kernel %;$\bL(\sN)$-valued function
\[
%\fL(z,\xi)=
\frac{(1-z^2)Q(z)-(1-\bar\xi^2)Q^*(\xi)-
(1-z\bar\xi)Q^*(\xi)(F-F^*)Q(z)-(z-\bar\xi) I}{z-\bar\xi}
%\quad |z|,|\xi|>1.
\]
is  nonnegative;
\item[(S4)] there exist a complex number $z_0$, $|z_0 |>1$, and a vector
$f\in\sN$, such that
\[
\frac{Q(z_0)-Q^*(z_0)-Q^*(z_0)(F-F^*)Q(z_0)}{z_0-\bar{z_0}}\ne
Q^*(z_0)Q(z_0)f.
\]
\end{enumerate}

If $T$ is a $qsc$-operator in the Hilbert space $\sH$, $\sN$ is a
subspace of $\sH$ such that $\ran (T-T^*)\subset\sN$, and $Q_T(z)$
is its $Q$-function defined by \eqref{QF}, then  the function
$Q_T(z)$ belongs to the class $\bQ(\sN)$, see \cite{AHS1}. The
converse statement is also true.

\begin{theorem} [\cite{AHS1}]\label{IP}
Let $Q(z)$ belong to $\bQ(\sN)$. Then there exist  Hilbert spaces
$\cH\supset\sN$, $\sN \ne \cH$, and a $\sN$-minimal $qsc$-operator
$T$ in $\sH$, such that  $\sN\supset \ker(T-T^*)$, and
\[
Q(z)=P_\sN(T-zI)^{-1}\uphar\sN, \quad  |z|>1.
\]
\end{theorem}

The next proposition gives connections between the transfer function
of a $pqs$-system  $\tau$ and the $Q$-function of the corresponding
$qsc$-operator $T$.

\begin{proposition}\label{QF&TR}
Let $\tau=\{T;\sN,\sN, \sH\}$ be a $pqs$-system.  Then the transfer
function $\Theta(\lambda)$ of $\tau$ and the $Q$-function of the
$qsc$-operator $T$  are connected by the following relations
\begin{equation}\label{QTheta}
Q(z)=\left(\Theta\left(\frac{1}{z}\right)-z I_\sN\right)^{-1}, \,\,
|z|>1;\quad \Theta(\lambda)=\frac{1
}{\lambda}I+Q^{-1}\left(\frac{1}{\lambda}\right), \,\, \lambda \in
\dD.
\end{equation}
\end{proposition}

\begin{proof}
With $W(z)=\Theta(1/z) -zI$, $|z|>1$, the resolvent $(T-z I)^{-1}$
  has the form
\[
\begin{pmatrix}  W^{-1}(z)&-W^{-1}(z)C(A-z I )^{-1}
\cr -(A-z I )^{-1}BW^{-1}(z) & (A-z I )^{-1}\left(I +BW^{-1}(z)C(A-z
I )^{-1}\right)\end{pmatrix}.
\]
Here the Schur-Frobenius formula has been applied; cf. \cite[Section
2.4]{AHS1}. It follows that
\[
P_\sN(T-zI)^{-1}\uphar \sN= W^{-1}(z).
\]
Thus, the relations \eqref{QTheta} hold.
\end{proof}

Observe that Theorem  \ref{recsys} and Theorem \ref{IP} imply that
$\Theta(\lambda)$ belongs to ${\bf S}^{qs}(\sN)$ if and only if
\[
Q(z)=\left(\Theta\left(\frac{1}{z}\right)-z I_\sN\right)^{-1}, \quad
|z|>1,
\]
belongs to $\bQ(\sN)$.

\subsection{Scalar functions of the class ${\bf
S}^{qs}$ and Jacobi matrices}

Let $l_2(\dN)$ and $l_2(\dN_0)$ be the Hilbert spaces of square
summable complex-valued sequences
\[
x=\{x_1,x_2,\ldots,x_k,\ldots\}, \quad
x=\{x_0,x_1,\ldots,x_k,\ldots\},
\]
considered as semi-infinite vector-columns, with the inner product
given by
\[
(x,y)=\sum\limits_{k=1}^\infty x_k\bar{y}_k, \quad
(x,y)=\sum\limits_{k=0}^\infty x_k\bar{y}_k,
\]
respectively. Clearly $\dC\oplus l_2(\dN)= l_2(\dN_0)$. Define the
vectors $\{\delta_k\}$, $k \in \dN_0$, by
\[
\delta_0=(1, 0, 0, \ldots)^T,\quad \delta_k = (0,\ldots, 0 , 1, 0,
0, \ldots)^T,\quad k \in \dN,
\]
so that $1$ is the $(k+1)$-st entry. Then the vectors $\{\delta_k\}$
form an orthonormal basis in  $l_2(\dN_0)$.

\begin{theorem}\label{Jacob}
Let the scalar function $\Theta(\lambda)$ belong to the class ${\bf
S}^{qs}$.
\begin{enumerate}\def\labelenumi{\rm (\roman{enumi})}
\item
If $\Theta(\lambda)$ is rational  with $n$ poles then any minimal
$pqs$-system $\tau=\{T,\dC,\dC,\sH\}$ with the transfer function
$\theta(\lambda)$ is unitarily equivalent to the $pqs$-system
\[
\tau_0=\{T_0,\dC,\dC,\dC^n\},
\]
where the operator $T_0$ in the Hilbert space $\dC\oplus \dC^n=
\dC^{n+1}$ with respect to the canonical basis
$\{\delta_k\}_{k=0}^{n}$ is given by the three-diagonal Jacobi
matrix
\begin{equation}
\label{1} T_0=\begin{pmatrix} \Theta(0) & a_0 & 0    & 0 & \cdot &
\cdot & \cdot \\
a_0 & b_1 & a_1 & 0 & \cdot &
\cdot & \cdot \\
0    & a_1 & b_2 & a_2 & \cdot &
\cdot & \cdot \\
\cdot & \cdot & \cdot & \cdot &
\cdot & \cdot & \cdot \\
\cdot & \cdot & \cdot & \cdot &
\cdot & \cdot & a_{n-1} \\
\cdot & \cdot & \cdot & \cdot & 0 & a_{n-1} & b_n
\end{pmatrix}.
\end{equation}

\item If $\Theta(\lambda)$ is not rational then any minimal
$pqs$-system $\tau=\{T,\dC,\dC,\sH\}$ with the transfer function
$\Theta(\lambda)$ is unitarily equivalent to the $pqs$-system
$$\tau_0=\{T_0,\dC,\dC,l_2(\dN)\},$$ where the operator $T_0$ in the
Hilbert space $\dC\oplus l_2(\dN)= l_2(\dN_0)$ with respect to the
canonical basis $\{\delta_k\}_{k=0}^\infty$ is given by the
semi-infinite three-diagonal Jacobi matrix
\begin{equation} \label{01}
  T_0=\begin{pmatrix} \Theta(0) & a_0 & 0 &0   & 0 &
\cdot &
\cdot  \\
a_0 & b_1 & a_1 & 0 &0& \cdot &
\cdot  \\
0    & a_1 & b_2 & a_2 &0& \cdot &
\cdot   \\
\cdot & \cdot & \cdot & \cdot & \cdot & \cdot & \cdot
\end{pmatrix}.
\end{equation}
\end{enumerate}
In both cases  $a_k>0$  and $b_k\in \dR$  for all relevant $k$, and
these numbers are uniquely determined by the function
$\Theta(\lambda)$.
\end{theorem}

\begin{proof}
Let $\Theta(\lambda) \in {\bf S}^{qs} $. By Proposition \ref{QF&TR}
the function $Q(z)=(\Theta(1/z) -z)^{-1}$, $|z|>1$, belongs to the
class ${\bf Q}$. Without loss of generality it can be assumed that
$\IM\Theta(0)\ge 0$. Since the function
$w(\lambda)=\Theta(\lambda)-\Theta(0)$ belongs to the
Herglotz-Nevanlinna class, the function $Q(z)$ has a holomorphic
continuation to the lower half-plane. If $\Theta(\lambda)$ is
rational with $n$ poles $\{\mu_k\}$, then these poles are simple and
belong to $(-\infty,-1)\cup(1,\infty)$. It follows that
\[
\Theta\left(\cfrac{1}{z}\right)=a\,\,\frac{\prod\limits_{k=1}^n(z-\lambda_k)}
{\prod\limits_{k=1}^n(z-\mu_k)},\quad \IM a\ge 0,
\]
and $Q(z)$ is rational with poles of total multiplicity $n+1$, in
the upper half-plane.

Let $\tau=\{T,\dC,\dC,\sH\}$ be a minimal $pqs$-system with the
transfer function $\Theta(\lambda)$. Then $Q(z)$ is the $Q$-function
of the corresponding $qsc$-operator $T,$ i.e.,
\[
Q(z)=\left((T-zI)^{-1}\bar 1,\bar 1\right), \quad |z|>1,
\]
where
\[
\bar 1=\begin{pmatrix}1\cr {\bf 0}\end{pmatrix}, \quad 1 \in \dC,
\quad {\bf 0} \in \sH \quad \mbox{(null-vector)}.
\]
The operator $T$ in the Hilbert space $\cH=\dC\oplus \sH$ is
dissipative ($\IM T(f,f)\ge 0$ for all $f\in\cH$), $\dim
\ran(T-T^*)\le 1$, and by Proposition~\ref{mpss} the operator $T$ is
$\dC$-minimal, i.e., $\cspan\{T^n \bar 1:\,n\in\dN_0\}=\cH$. If $\IM
\Theta(0)=0$, then $T$ is selfadjoint with the cyclic vector $\bar
1$. If $\IM \Theta(0)\ne 0$, then $T$ is a prime dissipative
operator with a rank-one imaginary part and, moreover, $\ran
(T-T^*)=\dC\oplus\{0\}$.

Note that $T$ is unitarily equivalent  to $T_0$ given by \eqref{1}
(in the case of a rational $Q(z)$) or by \eqref{01} (in the
opposite case) with respect the orthonormal basis $\{\delta_k\}$,
i.e., there exists a unitary operator $U\in{\bf L}(\cH,\dC^{n+1})$
or $U\in {\bf L}\left(\cH,l_2(\dN_0)\right)$ such that $UT=T_0 U$;
cf. \cite[Theorem 2.10, Theorem 6.1, Remark 6.2]{ArlTsek}.
Moreover,
\begin{equation}
\label{U} U\bar 1=\delta_0.
\end{equation}
Thus, for $|z|>1$ it follows that
$Q(z)=\left((T_0-zI)^{-1}\delta_0,\delta_0\right)$. The entries of
the matrix $T_0$  can be found by the continued-fraction expansion
of the function $Q(z)$:
\[
Q(z)= \frac{-1}{z-\Theta(0)}\;\raisebox{-3mm}{{\rm
+}}\;\frac{-a^2_0}{z-b_1}\;\raisebox{-3mm}{{\rm
+}}\;\frac{-a^2_1}{z-b_2} \;\raisebox{-3mm}{{\rm +}\,\ldots}\;
\raisebox{-3mm}{{\rm +}}\;\frac{-a^2_{n-1}}
{z-b_n}\;\raisebox{-3mm}{{\rm +}\,\dots},
\]
Note that $T_0$ is selfadjoint with the cyclic vector $\delta_0$
if $\IM\Theta(0)=0$, and that $T_0$ is a prime dissipative operator
with a rank-one imaginary part and $\ran(T_0-T_0^*)=\lin\{\delta_0\}$  if
$\IM\Theta(0)\ne 0$. From \eqref{01} it follows that
\[
T_0=\begin{pmatrix}D_0& B^*_0\cr
B_0&A_0\end{pmatrix}:\begin{pmatrix}\dC\\ \dC^n\end{pmatrix}\to
\begin{pmatrix}\dC\\ \dC^n\end{pmatrix}\,\, \mbox{or} \,\,
 T_0=\begin{pmatrix}D_0& B^*_0\cr
B_0&A_0\end{pmatrix}:\begin{pmatrix}\dC\\ l_2(\dN)\end{pmatrix}\to
\begin{pmatrix}\dC\\ l_2(\dN)\end{pmatrix},
\]
respectively, where $D_0=\Theta(0)$,
\[
B_0 1=\begin{pmatrix}a_0\cr 0\cr 0\cr\vdots\end{pmatrix}, \quad
B^*_0\begin{pmatrix}x_1\cr x_2\cr x_3 \cr \vdots\end{pmatrix}=a_0
x_1 \delta_0,
\]
and
\[
A_0\begin{pmatrix}x_1\cr x_2\cr\vdots\end{pmatrix}=\begin{pmatrix}
 b_1 & a_1 & 0 &0& \cdot &
\cdot  \\
 a_1 & b_2 & a_2 &0& \cdot & \cdot   \\
 \cdot & \cdot & \cdot & \cdot & \cdot & \cdot
\end{pmatrix}\begin{pmatrix}x_1\cr x_2\cr\vdots\end{pmatrix}.
\]
Decompose $T$ according to $\cH=\dC\oplus\sH$:
\[
T=\begin{pmatrix}D& B^*\cr B&A\end{pmatrix}:\begin{pmatrix}\dC\\
\sH\end{pmatrix} \to\begin{pmatrix}\dC\\ \sH\end{pmatrix}.
\]
Because \eqref{U} holds, the unitary operator $U$ takes the
following block operator matrix form
\[
U=\begin{pmatrix}1&0\cr 0& V\end{pmatrix}\begin{pmatrix}\dC\\
\sH\end{pmatrix} \to\begin{pmatrix}\dC\\
\dC^n\end{pmatrix}\quad \mbox{or} \quad  U=\begin{pmatrix}1&0\cr 0&
V\end{pmatrix}\begin{pmatrix}\dC\\ \sH\end{pmatrix} \to\begin{pmatrix}\dC\\
\l_2(\dN)\end{pmatrix},
\]
respectively. Hence, $VA=A_0V,$ $VB=B_0$, i.e., the $pqs$-systems
$\tau$ and $\tau_0$ are unitarily equivalent.
\end{proof}

Here is a simple example to illustrate the situation.

\begin{example}
Consider the scalar-valued function $\Theta(\lambda)$ defined  by
\[
\Theta(\lambda)=d+\frac{1-\sqrt{1-\lambda^2}}{2\lambda},\quad\lambda\in
\Ext\{(-\infty,-1]\cup[1,+\infty)\}.
\]
Then $\Theta(0)=d$ and $W(\lambda)=\Theta(\lambda)-\Theta(0)$ is
given by
\[
W(\lambda)=\frac{1-\sqrt{1-\lambda^2}}{2\lambda}, \quad\lambda\in
\Ext\{(-\infty,-1]\cup[1,+\infty)\}.
\]
Clearly $W(\lambda)$ is a Herglotz-Nevanlinna function. It follows that
\[
-2W\left(\frac{1}{z}\right)%%=-\frac{1}{W\left(\frac{1}{z}\right)}
=\sqrt{z^2-1}-z=\frac{1}{\pi}\,
\int\limits_{-1}^1\frac{\sqrt{1-t^2}\,dt}{t-z}, \quad
z\in\Ext[-1,1].
\]
Hence $W(1)=-W(-1)=1/2$. Moreover,  $\Theta(\lambda)$ belongs to the
class ${\bf S}^{qs}$ if and only if $|d|\le 1/2$; see \eqref{BW0}.
Assume that this
condition is satisfied. Consider the weighted Hilbert space
$L_2\left([-1,1],\rho(t)\right)$ with the weight function
\[
\rho(t)=\frac{2}{\pi}\, \sqrt{1-t^2}, \quad t \in [-1,1].
\]
Define the operator $A$ in $L_2\left([-1,1],\rho(t)\right)$ by
\[
(Af)(t)=tf(t),\quad f(t)\in L_2\left([-1,1],\rho(t)\right).
\]
Then $A$ is a selfadjoint contraction. The function $e_0(t)=1$,
$t\in [-1,1]$ belongs to $L_2\left([-1,1],\rho(t)\right)$ and
$\|e_0\|=1$.  Define the operator $B:\dC\to
L_2\left([-1,1],\rho(t)\right)$ by
\[
B c=\frac{1}{2}\,c\,e_0(t),\quad c\in\dC.
\]
Then
\[
B^*f(t)=\frac{1}{\pi}\,\int\limits_{-1}^1
f(t){\sqrt{1-t^2}\,dt},\quad f(t)\in L_2\left([-1,1],\rho(t)\right).
\]
Let $D$ be the multiplication by $d$ in the space $\dC$. One can
check that
\[
T=\begin{pmatrix}D&B^*\cr B&A\end{pmatrix}:\begin{pmatrix}\dC\\
L_2\left([-1,1],\rho(t)\right)\end{pmatrix} \to
\begin{pmatrix}\dC\\ L_2\left([-1,1],\rho(t)\right)\end{pmatrix}
\]
is a $qsc$-operator. Moreover, the corresponding $pqs$-system
\[
\tau=\left\{T;\dC,\dC,L_2\left([-1,1],\rho(t)\right)\right\}
\]
is minimal  and has the transfer function $\Theta(\lambda)$.  Since
the operator $A$ is unitarily equivalent to the Jacobi matrix (see
\cite{Ber})
\[
A_0 =\begin{pmatrix}0 & 1/2 & 0 &0& \cdot &
\cdot & \cdot \\
1/2    & 0 & 1/2 & 0 &0& \cdot &
\cdot & \cdot \\
\cdot & \cdot & \cdot & \cdot &
\cdot & \cdot & \cdot &\cdot \\
\end{pmatrix},
\]
the unitarily equivalent three-diagonal minimal $pqs$-system is of
the form
\[
 T_0=\begin{pmatrix} d & 1/2 & 0  &0  & 0 & \cdot &
\cdot & \cdot \\
1/2 &0 & 1/2 & 0 &0& \cdot &
\cdot & \cdot \\
0    & 1/2 & 0 & 1/2 &0& \cdot &
\cdot & \cdot \\
\cdot & \cdot & \cdot & \cdot &
\cdot & \cdot & \cdot &\cdot \\
\end{pmatrix}.
\]
\end{example}

\end{document}